\documentclass[12pt,draft]{amsart}
\usepackage[T1]{fontenc}
\usepackage{times}

\makeatletter

\newcommand{\LyX}{L\kern-.1667em\lower.25em\hbox{Y}\kern-.125emX\spacefactor1000}

\theoremstyle{plain}    
\newtheorem{thm}{Theorem}[section]
\numberwithin{equation}{section} 
\numberwithin{figure}{section} 
\theoremstyle{plain}    
\newtheorem*{thm*}{Theorem} 
\theoremstyle{plain}    
\newtheorem{cor}[thm]{Corollary} 
\theoremstyle{plain}    
\newtheorem{lem}[thm]{Lemma} 
\theoremstyle{plain}    
\newtheorem{prop}[thm]{Proposition} 
\theoremstyle{definition}

\theoremstyle{remark}
\newtheorem{rem}[thm]{Remark}
\theoremstyle{remark}
\newtheorem{rems}[thm]{Remarks}
\theoremstyle{remark}    
\newtheorem{claim}[thm]{Claim}
\theoremstyle{remark}    
\newtheorem{question}[thm]{Question}
\theoremstyle{remark}    
\newtheorem*{acknowledgement*}{Acknowledgement} 

\theoremstyle{plain}    
\newtheorem{subthm}{Theorem}[subsection]
\theoremstyle{plain}    
\newtheorem{subcor}[subthm]{Corollary} 
\theoremstyle{plain}    
\newtheorem{sublem}[subthm]{Lemma} 
\theoremstyle{plain}    
\newtheorem{subprop}[subthm]{Proposition} 
\theoremstyle{definition}
\newtheorem{subdefn}[subthm]{Definition}
\theoremstyle{remark}
\newtheorem{subrem}[subthm]{Remark}
\theoremstyle{remark}    

\theoremstyle{remark}    
\newtheorem{subquestion}[subthm]{Question}
\theoremstyle{plain}    

\newcommand{\Id}{\operatorname{Id}}
\newcommand{\id}{\operatorname{id}}
\newcommand{\Ad}{\operatorname{Ad}}
\newcommand{\HT}{\operatorname{ht}}
\newcommand{\Span}{\operatorname{span}}
\newcommand{\Aut}{\operatorname{Aut}}

\newcommand{\rcp}{\operatorname{rcp}}
\newcommand{\rank}{\operatorname{rank}}

\usepackage{amssymb}
\usepackage{amscd}

\makeatother

\begin{document}

\title{topological entropy of free product automorphisms}

\date{October 10, 2000}

\author{N.P. Brown, K. Dykema and D. Shlyakhtenko}

\address{Department of Mathematics, University of California, Berkeley, CA 94720}

\email{nbrown@math.berkeley.edu}

\address{Department of Mathematics, Texas A\&M University, College Station, TX 77843}

\email{Ken.Dykema@math.tamu.edu}

\address{Department of Mathematics, University of California, Los Angeles, CA 90095}

\email{shlyakht@math.ucla.edu}
\thanks{N.B. and D.S. partially supported by NSF postdoctoral fellowships.
K.D. partially supported by NSF Grant No.\ DMS 0070558. N.B. is currently an MSRI Postdoctoral Fellow.}

\begin{abstract}
Using free probability constructions involving the Cuntz--Pimsner \linebreak[4]
\(C^{*} \)--algebra we show that the topological entropy \( \HT (\alpha
*\beta ) \) of the free product of two automorphisms is given by the
maximum \( \max (\HT (\alpha ),\HT (\beta )) \).  As applications, we
show in full generality that free shifts have topological entropy
zero. We show that any separable nuclear \( C^{*} \)--dynamical system
can be covariantly embedded into \( \mathcal{O}_{2} \) and \(
\mathcal{O}_{\infty } \) in an entropy--preserving way. It follows
that any nuclear simple purely infinite \( C^{*} \)--algebra admits an
automorphism with any given value of entropy. We also show that the free
product of two automorphisms satisfies the Connes--Narnhofer--Thirring
variational principle, if the two automorphisms do. 
\end{abstract}
\maketitle

\section{Introduction.}

The study of non--commutative topological dynamical systems, in other
words, automorphisms of \( C^{*} \)--algebras, goes back a long
way. Among the most important invariants of such systems is entropy,
first introduced in the operator algebra context by Connes and
St\o{}rmer in \cite{connes-stormer:entropy}. (See
\cite{stormer:entropysurvey} for an extensive survey of this subject.)
In this paper, we deal with the non--commutative generalization of
topological entropy, discovered by Voiculescu \cite{dvv:topentropy}
for nuclear \( C^{*} \)--algebras and extended by Brown
\cite{brown:entropy} to the exact case. Although this entropy has been
computed in many examples, we are very far from a good understanding
of its behavior. For example, to this day we do not know in full
generality the exact value of the topological entropy \( \HT (\alpha
\otimes \beta ) \) of a tensor product of two automorphisms.  In
general the inequality \( \HT (\alpha \otimes \beta )\le\HT (\alpha
)+\HT (\beta ) \) is known.  Equality is only known to hold when both
\( \alpha \) and \( \beta \) satisfy a CNT--variational principle.

In this paper, we address four general types of questions regarding topological
entropy.

The first question concerns the behavior of topological entropy under
free products.  The precise question is the value of \( \HT (\alpha
*\beta ) \) of the topological entropy of the (possibly amalgamated)
reduced free product automorphism \( \alpha *\beta \).  The main
theorem of this paper states that \( \HT (\alpha *\beta )=\max (\HT
(\alpha ),\HT (\beta )) \), if the free product is with amalgamation
over a finite dimensional \( C^{*} \)--algebra.  One surprising feature
of this result is that the answer is precise --- this is to be
contrasted with the situation for tensor products. Although the free
product of \( C^{*} \)--algebras is more complex than the tensor
product, it seems to behave more like the direct sum for the purposes
of entropy, which is, curiously, close to its behaviour for the purposes
of  K--theory. One direct consequence of our result is that the free
shift on an arbitrary reduced free product \( *_{i\in \mathbb
{Z}}(A,\phi ) \) has zero entropy.

Giving a bound on the topological entropy of a free product of two
automorphisms can be considered a refinement of the
result~\cite{dykema:exact} that the reduced amalgamated free product
of two exact \( C^{*} \)--algebras is exact. Indeed, exactness of a
\( C^{*} \)--algebra \( B \) is the statement that \( \HT (\id :B\to B)<+\infty  \).
The reader will find that our proof is related to the argument given
in
\cite{shlyakht-dykema:exactness} for the exactness of reduced amalgamated free
products. 

The second general question concerns the Connes--Narnhofer--Thirring
(CNT) variational principle. The Connes--Narnhofer--Thirring (CNT)
entropy \cite{connes-narnhofer-thirring:entropy} is a generalization
to  non--commutative measure spaces of the classical
Kolmogorov--Sinai entropy. An automorphism \( \alpha \) is said to
satisfy the CNT variational principle if its topological entropy is
equal to the supremum of the CNT entropy computed with respect to all
invariant states. Although this principle fails for general
non--commutative dynamical systems, we show that if two automorphisms
\( \alpha \in \Aut (A) \) and \( \beta \in \Aut (B) \) satisfy the 
CNT variational principle, then so does \( \alpha *\beta \).

The third general question concerns embeddings of dynamical
systems. Kirchberg has shown that any separable exact \( C^{*}
\)--algebra can be embedded into the Cuntz algebra on two
generators. We show that any nuclear \( C^{*} \)--dynamical system can
be covariantly embedded into the Cuntz algebra  \( \mathcal{O}_\infty \)  in an
entropy--preserving way. 

The last question concerns the possible values of entropy that can be
attained by automorphisms of a given \( C^{*} \)--algebra \( A \),
i.e., the set \( TE(A)=\HT (\Aut (A)) \).  This is clearly an
invariant of \( A \). For instance, \( TE(A)=\{+\infty \} \) is \( A
\) is not exact; \( TE(A)=\{0\} \) if \( A \) is finite dimensional.
We show that \( TE(\mathcal{O}_{\infty })=[0,+\infty ] \).  
Thanks to Kirchberg's absorption results, this
implies that any separable purely infinite nuclear simple \( C^{*} \)
--algebra \( A \) admits an automorphism with any given value of
entropy; i.e., \( TE(A)=[0,+\infty ] \).

The main result of this paper, computing entropy of amalgamated free products
of automorphisms, applies only to the case that amalgamation takes place over
a finite dimensional subalgebra. Our results are likely to extend to the
case of amalgamation over an arbitrary \( C^{*} \)--subalgebra, if the following
question can be answered in the affirmative:

\begin{question}Let \( 0\hookrightarrow I\to E\stackrel{\pi }{\to }B\to 0 \)
be an exact sequence of \( C^{*} \)--algebras, which is split (i.e.,
there exists a \( * \)--homomorphism \( \eta :B\to E \), so that \( \pi
\circ \eta =\id _{B} \) ). Let \( \alpha \) be an automorphism of \( E
\), so that \( \alpha (I) = I \); denote the resulting
automorphism of \( B \) by \( \hat{\alpha } \). Assume that \( \eta
\circ \hat{\alpha }=\alpha \circ \eta \).  Is it true that
\( \HT (\alpha )\leq \max (\HT (\hat{\alpha }),\HT(\alpha |_{I})) \)?\end{question}

\begin{acknowledgement*}
An early part of this work was carried out at Institute Henri Poincar\'e, Paris,
France, to which the authors would like to express their warmest gratitude for
the exciting and encouraging atmosphere. 
We would like to thank also N.\ Ozawa
and E.\ Germain for helpful conversations. 
K.D.\ would like to thank Uffe Haagerup and Odense University for their kind hospitality
when part of this work was carried out.
\end{acknowledgement*}

\section{Preliminaries.}

\subsection{A review of topological entropy.}

Let \( A \) be an exact \( C^{*} \) algebra and let \( \omega \subset A \)
be a finite subset. Fix a faithful representation \( \pi :A\to B(H) \). Define
the set \( CPA(A) \) to be the set of triples \( (\Phi ,\Psi ,X) \), where
\( X \) is a finite dimensional \( C^{*} \)--algebra, and \( \Phi :A\to X \),
\( \Psi :X\to B(H) \) are contractive completely--positive maps. For \( \epsilon >0 \),
define
\[
\rcp (\omega ,\epsilon )=
\inf (\rank (X):(\Phi ,\Psi ,X)\in CPA(A),\quad \forall a\in \omega\;
\Vert \Psi \circ \Phi (a)-\pi (a)\Vert <\epsilon  ).\]
(This quantity is independent of the choice of \( \pi  \), cf. \cite{brown:entropy}).
If \( \alpha \in \Aut (A) \) is an automorphism, then its topological entropy
\( \HT (\alpha ) \) is defined as
\[
\sup _{\omega \subset A}\sup _{\epsilon }\limsup _{n\to \infty }\frac{1}{n}\log \rcp (\omega \cup \cdots \cup \alpha ^{n-1}(\omega ),\epsilon ).\]
This definition, obtained in \cite{brown:entropy}, gives the same quantity
as Voiculescu's original definition of dynamical topological entropy \cite{dvv:topentropy}
for nuclear \( C^{*} \)--algebras \( A \). We summarize below some properties
of \( \HT (\alpha ) \), which we will need in this paper. The proofs can be
found in \cite{brown:entropy}.

\begin{thm*}
Let \( \alpha \in \Aut (A) \) be as above. Then
\begin{enumerate}
\item \( \HT  \) is monotone: if \( B\subset A \), and \( \alpha (B)\subset B \),
then \( \HT (\alpha |_{B})\leq \HT (\alpha ) \).
\item If \( A=\overline{\bigcup A_{n}} \), with \( A_{n} \subset A_{n+1} \) subalgebras and \( \alpha (A_{n})\subset A_{n} \),
then \( \HT (\alpha )=\lim _{n}\HT (\alpha |_{A_{n}}) \).
\item If \( \beta \in \Aut (B) \), let \( \alpha \otimes \beta \in \Aut (A\otimes_{\min} B) \)
be the tensor product automorphism. Then \( \HT (\alpha \otimes \beta
)\leq \HT (\alpha )+\HT (\beta ) \).  If \( A \) contains an \( \alpha \)--invariant 
projection, then \( \HT (\alpha \otimes \beta )\geq \HT (\beta ) \).
\item If \( \beta \in \Aut (B) \), let \( \alpha \oplus \beta \in \Aut (A\oplus B) \)
be the direct sum automorphism. Then \( \HT (\alpha \oplus \beta )=\max (\HT (\alpha ),\HT (\beta )) \).
\item If \( \beta \in \Aut (A) \) commutes with \( \alpha  \), i.e., \( \beta \circ \alpha =\alpha \circ \beta  \),
then \( \alpha  \) extends to the obvious automorphism \( \bar{\alpha }\in \Aut (A\rtimes _{\beta }\mathbb {Z}) \).
Then \( \HT (\bar{\alpha })=\HT (\alpha ) \).
\end{enumerate}
\end{thm*}

\subsection{Amalgamated free products.}
\label{subsec:fp}

Let \( D \) be a unital \( C^{*} \)--algebra. Recall that a \( D \)--\emph{valued
non--commutative probability space} is a \( C^{*} \)--algebra \( A \), containing
\( D \) as a unital subalgebra, and endowed with a conditional expectation
\( E:A\to D \). 

Let now \( A_{1} \) and \( A_{2} \) be two unital \( C^{*} \)--algebras with
a common unital subalgebra \( D \) and with conditional expectations \( E_{j}:A_{j}\to D \),
so that the GNS representation associated to each \( E_{j} \) is faithful.
Then there exists a \( C^{*} \)--algebra \( A \), generated by \( A_{1} \)
and \( A_{2} \) (in such a way that the copies of \( D \) inside \( A_{j}\subset A \)
are identified), a conditional expectation \( E:A\to D \), and such that \( A_{j}\subset A \)
are \emph{free with respect to \( E \)} (see \cite{DVV:book}), and
\( E \) gives rise to a faithful GNS representation of \( A \). The \( C^{*} \)--algebra
\( A \) is the \emph{reduced amalgamated free} product over \( D \) of \( (A_{1},E_{1}) \)
and \( (A_{2},E_{2}) \) and is denoted
\begin{equation}
\label{eq:amalgfp}
(A,E)=(A_1,E_1)*_D(A_2,E_2).
\end{equation}
When the conditional expectations $E_j$ are clear from context we may write $A=A_1*_DA_2$.
In the case that \( D=\mathbb {C}1\subset A_{1},A_{2} \),
then $E_1$, $E_2$ and $E$ are states and \( A \) is the reduced free product of \( A_{1} \)
and \( A_{2} \), denoted simply
\begin{equation}
\label{eq:fp}
(A,E)=(A_1,E_1)*(A_2,E_2).
\end{equation}
Moreover, in the situation of a free product~\eqref{eq:amalgfp} or~\eqref{eq:fp}, one has
conditional expectations $\Phi_1:A\to A_1$ and $\Phi_2:A\to A_2$ such that $E_j\circ\Phi_j=E$.

\begin{subdefn}
Let \( (A,E:A\to D) \) be a \( D \)--valued non--commutative probability space.
We say that \( \alpha \in \Aut (A) \) is an automorphism of this space, if
(i) \( \alpha (D)=D \) and (ii) \( \alpha \circ E=E\circ \alpha  \).
\end{subdefn}
Notice that in the case that \( D=\mathbb {C} \), an automorphism of a \( D \)--probability
space \( A \) is simply an automorphism of \( A \), fixing the state \( E:A\to D=\mathbb {C} \).

Assume now that \( (A_{1},\;E:A_{1}\to D) \) and \( (A_{2},\;E_{2}:A_{2}\to D) \)
are \( D \)--probability spaces. Assume further that \( \alpha _{j}\in \Aut (A_{j}) \)
are \( D \)--space automorphisms, and \( \alpha _{1}|_{D}=\alpha _{2}|_{D} \).
Then is  a unique automorphism
\( \alpha _{1}*\alpha _{2}:A_1*_DA_2\to A_1*_DA_2 \), satisfying
\( (\alpha _{1}*\alpha _{2})|_{A_{j}\subset A_{1}*_{D}A_{2}}=\alpha _{j} \),
\( j=1,2 \).

We will need the following theorem of Blanchard and Dykema.
\begin{thm*}
\cite{blanch-dyk}
Let $D$ and $\widetilde D$ be unital $C^*$--algebras, let
$(A_1,E_1)$, $(A_2,E_2)$ be $D$--probability spaces and let
$(\widetilde{A}_1,\widetilde{E}_1)$ and $(\widetilde{A}_2,\widetilde{E}_2)$
be $\widetilde D$--probability spaces.
Assume that the GNS representations associated to $E_j$ and $\widetilde E_j$ are faithful ($j=1,2)$.
Let
\begin{align*}
(A,E)&=(A_1,E_1)*_D(A_2,E_2) \\
(\widetilde A,\widetilde E)&=(\widetilde A_1,\widetilde E_1)*_{\widetilde D}(\widetilde A_2,\widetilde E_2).
\end{align*}
Suppose that $\pi_D:D\to\widetilde D$
is a (not necessarily unital) injective $*$--homomorphism
and that there are injective $*$--homomorphisms
$\pi_1:A_1\to\widetilde{A}_1$ and $\pi_2:A_2\to\widetilde{A}_2$
such that $\pi_D\circ E_j=\widetilde E_j\circ\pi_j$.
Then there is an injective $*$--homomorphism $\pi:A\to\widetilde A$ such that $\pi|_{A_j}=\pi_j$ ($j=1,2$) and
$\widetilde E\circ\pi=E$.
\end{thm*}

\subsection{\protect\( K\protect \)--theory for reduced free products.\label{sec:Ktheoryfreeprod}}

Thanks to fundamental work of E. Germain on K--theory of free products ,
one has the following six--term exact sequence for free products of nuclear \( C^{*} \)--algebras.

\begin{thm*}
\cite{germain:KKfullfreeprod}, \cite{germain:KKredfreeprod}, \cite{germain:Kcommutatorideal}
 Let \( A \) and \( B \) be unital nuclear \( C^{*} \)--algebras and \( \phi \in S(A),\psi \in S(B) \)
be states with faithful GNS representations.  Let \( C=(A,\phi )*(B,\psi ) \)
and \( \iota _{A}:A\hookrightarrow C \) and \( \iota _{B}:B\hookrightarrow C \)
denote the canonical inclusions.  Then there is an exact sequence $$\begin{CD}
\mathbb{Z} @>1\mapsto [1_A] \oplus -[1_B]>> K_0(A)\oplus K_0(B) @>(\iota_A)_*+(\iota_B)_*>>K_0(C)\\
@AAA @. @VVV\\
K_1(C) @<(\iota_A)_*+(\iota_B)_*<< K_1 (A)\oplus K_1 (B) @<<< 0.
\end{CD}$$ Moreover, if both \( A \) and \( B \) satisfy the Universal Coefficient
Theorem of Rosenberg and Schochet (see e.g. \cite{blackadar:book}) then so
does \( C \). 
\end{thm*}

\subsection{Cuntz--Pimsner algebras.}

\label{sec:CP-algebras}

Let \( A \) be a \( C^{*} \)--algebra, and  \( H \) be a Hilbert bimodule over $A$.
Assume that \( H \) is \emph{full}, i.e., \( \langle H,H\rangle _{A} \) is
dense in \( A \). Let
\[
F(H)=A\oplus \bigoplus _{n\geq 1}H^{\otimes _{A}n},\]
and for \( \xi \in H \), let
\begin{eqnarray*}
l(\xi ):F(H) & \to  & F(H),\\
l(\xi )\cdot \xi _{1}\otimes \cdots \otimes \xi _{n} & = & \xi \otimes \xi _{1}\otimes \cdots \otimes \xi _{n}.
\end{eqnarray*}
The Cuntz--Pimsner \( C^{*} \)--algebra \( E(H) \) (cf. \cite{pimsner}) is
then defined as \( C^{*}(l(\xi ):\xi \in H). \) It is not hard to see that
\[
l(\xi )^{*}l(\zeta )=\langle \xi ,\zeta \rangle _{A},\]
and hence \( A=\overline{\langle H,H\rangle }\subset E(H) \), acting on the left of $F(H)$.
The projection
from \( F(H) \) onto \( A\subset F(H) \) gives rise to a conditional expectation
\( E:E(H)\to A \). We summarize some of the properties of \( E(H) \) below;
the proofs can be found in \cite{pimsner}, \cite{speicher:thesis} and \cite{shlyakht-dykema:exactness}.

\begin{thm*}
Let \( H \) and \( A \) be as above. Then
\begin{enumerate}
\item If \( K\subset H \) is a Hilbert \( A \)--subbimodule which is full, then
$E(K)$ is canonically isomorphic to
\( C^{*}(l(\xi ):F(H)\to F(H)\, \textrm{s}.\textrm{t}.\, \xi \in K)\subset E(H) \).
\item Let \( H' \) be another full Hilbert \( A \)--bimodule, then \( E(H\oplus H')=E(H)*_{A}E(H') \),
where the reduced free product is taken with respect to the canonical conditional
expectations from \( E(H) \) (respectively, \( E(H') \)) onto \( A \).
\item Assume that \( \bar{\alpha }:H\to H \) is a linear map, so that for some \( \alpha \in \Aut (A) \),
\begin{eqnarray}
\bar{\alpha }(a_{1}\cdot \xi \cdot a_{2})=\alpha (a_{1})\cdot \bar{\alpha }(\xi )\cdot \alpha (a_{2}), &  & a_{1},a_{2}\in A,\, \xi \in H\label{eqn:Bogoljubovaxioms1} \\
\langle \bar{\alpha }(\xi _{1}),\bar{\alpha }(\xi _{2})\rangle _{A}=\alpha (\langle \xi _{1},\xi _{2}\rangle _{A}), &  & \xi _{1},\xi _{2}\in H.\label{eqn:Bogoljubovaxioms2} 
\end{eqnarray}
Then there is a unique automorphism \( E(\bar{\alpha }) \) of the \( A \)--probability
space \( (E(H),E:E(H)\to A) \), so that \( E(\bar{\alpha })|_{A\subset E(H)}=\alpha  \),
and \( E(\bar{\alpha })(l(\xi ))=l(\bar{\alpha }(\xi )) \). The automorphism
\( E(\bar{\alpha }) \) is called the Bogoljubov automorphism associated to
\( \bar{\alpha } \).
\end{enumerate}
\end{thm*}

\subsection{\label{sec:exampleBogoljubov}Some examples of Bogoljubov automorphisms.}

In the course of proving the main theorem of the paper, we will encounter a
particular class of Bogoljubov automorphisms of \( E(H) \), which we will presently
describe. Let \( D \) be a \( C^{*} \)--algebra and \( (A,E:A\to D) \) be
a \( D \)--probability space. Let
\[
K^{o}=A\otimes A\]
(algebraic tensor product), and endow \( K^{o} \) with the \( A \)--valued
inner product given by
\[
\langle a\otimes b,a'\otimes b'\rangle _{A}=b^{*}E(a^{*}a')b',\qquad a,a',b,b'\in A.\]
Denote by \( K_{D} \) the Hilbert \( A \)--bimodule obtained from \( K^{o} \)
after separation and completion.
Another description of $K_D$ is as the internal tensor product
$K_D=L^2(A,E)\otimes_DA$.
Notice that any automorphism \( \alpha  \)
of the \( D \)--probability space \( A \) extends to a linear map \( \bar{\alpha }:K_{D}\to K_{D} \),
satisfying equations (\ref{eqn:Bogoljubovaxioms1}) and (\ref{eqn:Bogoljubovaxioms2}). 

\begin{thm*}
\cite{shlyakht:amalg} \( E(K_{D})\cong (A,E)*_{D}(D\otimes \mathcal{T},\id \otimes \psi ) \),
where \( \mathcal{T} \) denotes the Toeplitz algebra (the algebra generated
by the unilateral shift on \( \ell ^{2}(\mathbb {N}) \)), and \( \psi :\mathcal{T}\to \mathbb {C} \)
is the vacuum state (corresponding to the vector \( \delta _{1}\in \ell ^{2}(\mathbb {N}) \)).
Moreover, \( E(\bar{\alpha }) \) corresponds in this isomorphism to \( \alpha *((\alpha |_{D})\otimes \id ) \).
\end{thm*}
Let \( (A,E:A\to D) \) be a \( D \)--probability space, and let \( \alpha  \)
be an automorphism of this \( D \)--probability space. Assume that \( D \)
is finite dimensional. Let \( \phi  \) be an \( \alpha |_{D} \)--invariant
faithful trace on \( D \). Consider the Hilbert \( A,A \)--bimodule
\[
H=L^{2}(A,\phi \circ E)\otimes _{\mathbb {C}}A,\]
together with the vector \( \xi =1\otimes 1\in H \) and the inner product
\[
\langle a_{1}\otimes a_{2},b_{1}\otimes b_{2}\rangle =\langle a_{1},b_{1}\rangle _{L^{2}(A,\phi \circ E)}a_{2}^{*}b_{2}.\]
We will henceforth also write $\phi$ to mean $\phi\circ E$.
Let \( U:H\to H \) be given by
\[
U(x\otimes a)=V(x)\otimes \alpha (a),\]
where \( V:L^{2}(A,\phi )\to L^{2}(A,\phi ) \) is the unitary induced by \( \alpha :A\to A \).
Then $UaU^*=\alpha(a)$ where $a$ and $\alpha(a)$ act on the left of $H$.

\begin{sublem}
\label{lemma:existenceexpectationvector}There exists a vector \( \zeta \in H \),
with the following properties:
\begin{enumerate}
\item \( d\zeta =\zeta d \) for all \( d\in D \) 
\item \( U\zeta =\zeta  \)
\item \( \langle \zeta ,a\zeta \rangle =E(a) \) for all \( a\in A \).
\end{enumerate}
\end{sublem}
\begin{proof}
Let \( U(D) \) denote the unitary group of \( D \), endowed with Haar measure
\( \mu  \). Let
\[
\zeta '=\int _{u\in U(D)}u\xi u^{*}d\mu (u).\]

For each \( w\in U(D) \), we get
\[
w\zeta '=w\int _{u\in U(D)}u\xi u^{*}d\mu (u)=\left( \int _{uw\in U(D)}wu\xi uw^{*}d\mu (u)\right) w=\zeta 'w,\]
so $d\zeta'=\zeta'd$ for all $d\in D$.
Furthermore,
\begin{eqnarray*}
U\zeta '=\int _{u\in U(D)}U(u\otimes u^{*})d\mu (u) & = & \int _{u\in U(D)}\alpha (u)\cdot 1\otimes \alpha (u^{*})d\mu (u)\\
 & = & \int _{u\in U(D)}\alpha (u)\xi \alpha (u^{*})d\mu (u)=\zeta '.
\end{eqnarray*}
Lastly, for \( a\in A \), set
\begin{eqnarray*}
\Phi (a)=\langle \zeta ',a\zeta '\rangle  & = & \int _{u,v\in U(D)}\langle u\otimes u^{*},av\otimes v^{*}\rangle d\mu (u)d\mu (v)\\
 & = & \int _{u,v}u\phi (u^{*}av)v^{*}d\mu (u)d\mu (v)\\
 & = & \int _{u,v}\phi (vu^{*}a)uv^{*}d\mu (u)d\mu (v)\\
 & = & \int _{w=uv^{*},v}\phi (wa)w^{*}d\mu (w)d\mu (v)\\
 & = & \int _{w}\phi (wa)w^{*}d\mu (w).
\end{eqnarray*}
Assume now that \( E(a)=0 \), i.e., \( \phi (da)=0 \) for all \( d\in D \).
Then
\[
\Phi (a)=\int _{w}\phi (wa)w^{*}d\mu (w)=0.\]
Since whenever \( v\in U(D) \),
\[
\Phi (v)=\int _{w}\phi (wv)w^{*}d\mu (w)=\int _{u=wv}\phi (u)vu^{*}d\mu (u)=v\Phi (1),\]
and any \( d\in D \) is a linear combination of unitaries from \( U(D) \),
we get that
\[
\Phi (d)=d\Phi (1),\qquad \forall d\in D.\]
Since for any \( a\in A \), \( a=(a-E(a))+E(a) \) and \( E(a)\in d \), we
get that
\[
\Phi (a)=\Phi (a-E(a))+\Phi (E(a))=E(a)\Phi (1).\]
We also have
\( d\Phi (1)=\Phi (d)=\langle \zeta ',d\zeta '\rangle =\langle \zeta ',\zeta 'd\rangle =\Phi (1)d \).
Hence \( \Phi (1) \) is in the center of \( D \); moreover, \( \Phi (1)=\langle \zeta ',\zeta '\rangle \geq 0 \)
and
\( \alpha (\Phi (1))=\alpha (\langle \zeta ',\zeta '\rangle )
=\langle U\zeta ',U\zeta '\rangle =\langle \zeta ',\zeta '\rangle =\Phi (1) \).
We claim that \( \Phi (1) \) is invertible. Writing \( D=\oplus M_{n_{k}} \)
as a direct sum of matrix algebras, we find that \( \Phi (1)=\int _{u\in U(D)}\phi (u)u^{*}d\mu (u)=\oplus \Phi _{k} \),
where \( \Phi _{k}=\alpha _{k}\cdot 1_{M_{n_{k}}}\cdot \int _{u\in U(n_{k})}\textrm{Tr}(u)u^{*}d\mu (u) \),
and \( \alpha _{k}\ge0 \) is related to the value of \( \phi  \) on the minimal
projection of \( D \) corresponding to the \( k \)-th matrix summand in the
direct sum decomposition of \( D \). To show that \( \Phi (1) \) is invertible,
is it sufficient to show that \( c_{k}=\int _{u\in U(k)}\textrm{Tr}(u)u^{*}d\mu (u) \)
is a strictly positive scalar for all \( k\geq 1 \). Repeating the argument
above with \( D \) replaced by \( M_{k} \), we find that \( c_{k} \) is in
the center of \( M_{k} \) and is non--negative. Furthermore,
\[
\textrm{Tr}(c_{k})=\int _{u\in U(k)}\textrm{Tr}(u)\textrm{Tr}(u^{*})d\mu (u)=\int _{u\in U(k)}|\textrm{Tr}(u)|^{2}d\mu (u)>0,\]
since the subset of \( u\in M_{k} \) with \( \textrm{Tr}(u)\neq 0 \) has non--zero
measure.

Now let \( \zeta =\Phi (1)^{-1/2}\zeta ' \). Then \( U\zeta =U\Phi (1)^{-1/2}U^{*}U\zeta '=\Phi (1)^{-1/2}\zeta '=\zeta  \);
for all \( d\in D \), \( d\zeta =d\Phi (1)^{-1/2}\zeta '=\Phi (1)^{-1/2}d\zeta '=\Phi (1)^{-1/2}\zeta 'd=\zeta d \);
and for all \( a\in A \),
\begin{eqnarray*}
\langle \zeta ,a\zeta \rangle = 
\langle \Phi (1)^{-1/2}\zeta ',a\Phi (1)^{-1/2}\zeta '\rangle  & = & 
\langle \zeta ',\Phi (1)^{-1/2}a\Phi (1)^{-1/2}\zeta '\rangle \\ 
& = & \Phi (\Phi (1)^{-1/2}a\Phi (1)^{-1/2})\\ 
& = & \Phi (1)^{-1/2}E(a)\Phi (1)^{-1/2}\Phi (1)\\ 
& = & E(a),
\end{eqnarray*}
as desired. 
\end{proof}
\begin{subcor}
\label{Corollary:embeddingToeplitz} Let \( \alpha  \) be an automorphism of
a \( D \)--probability space \( (A,E:A\to D) \) with \( \dim(D) < \infty\). Assume that \( E:A\to D \)
gives rise to a faithful GNS representation. Let \( \phi  \) be an \( \alpha  \)--invariant
trace on \( D \). Consider the automorphism \( \beta =\alpha *\id  \) on \( B=(A,\phi \circ E)*\mathcal{T} \),
where \( \mathcal{T} \) is the Toeplitz algebra with its vacuum state. Consider
the automorphism \( \gamma =\alpha *(\alpha |_{D}\otimes \id ) \) on \( C=(A,E)*_{D}(D\otimes \mathcal{T}) \).
Then there exists a covariant embedding of \( (C,\gamma ) \) into \( (B,\beta ) \).
\end{subcor}
\begin{proof}
Let \( H \) be as above. Consider the Cuntz--Pimsner \( C^{*} \)--algebra \( E(H) \),
generated by the operators \( l(h):h\in H \) and \( A \). Then \( B=E(H)=C^{*}(A,\, l(h):h\in H)=C^{*}(A,l(\xi )) \).
The automorphism \( \beta  \) is identified with the Bogoljubov automorphism
\( E(U) \):
indeed we have \( \beta (a)=\alpha (a) \), \( a\in A \) and
\( E(U)(l(a\otimes b))=l(U(a\otimes b))=l(\alpha (a)\otimes \alpha (b))=\alpha (a)l(\xi )\alpha (b)=\beta (al(\xi )b) \).
Choose \( \zeta  \) as in Lemma \ref{lemma:existenceexpectationvector}. Let
\( L=l(\zeta ) \). Then \( L^{*}aL=\langle \zeta ,a\zeta \rangle =E(a) \).
Let \( K=\overline{A\zeta A}\subset H \). It follows that \( C^{*}(A,L)=C^{*}(aLb:a,b\in A)=C^{*}(l(h):h\in K)=E(K) \).
It is easily seen that the map \( a\zeta b\mapsto a\otimes b \) defines an
isomorphism of \( K \) with the \( A \) bimodule \( K_{D} \) defined in \S\ref{sec:exampleBogoljubov}.
Hence by the results of that section, \( C^{*}(A,L)\cong C \). Moreover, since
\begin{eqnarray*}
E(U)(a)=\alpha (a),\quad a\in A &  & \\
E(U)(L)=l(U(\zeta ))=l(\zeta )=L &  & 
\end{eqnarray*}
we get that \( E(U)|_{C^{*}(A,L)}=\gamma  \).
\end{proof}

\section{Toplogical entropy in certain extensions}

The main result of this section gives an estimate of topological entropy in
certain extensions.  We begin with a lemma which gives a way of estimating the
\( \delta  \)--rank of finite sets in an extension.   

Consider a short exact sequence \( 0\to I\to E\stackrel{\pi }{\to }B\to 0 \)
 where \( E \) is a unital exact \( C^{*} \)--algebra and assume that there exists
a unital completely positive splitting \( \rho :B\to E \) (i.e. \( \pi \circ \rho =\id _{B} \)). 

\begin{lem}
\label{label:rcpestimateextension}Let \( I,E,B,\pi  \) and \( \rho  \) be
as above.  Given \( \epsilon >0 \) there exists a \( \delta >0 \) with the
following property: if \( \omega =\{\rho (b_{1})+x_{1},\ldots ,\rho (b_{s})+x_{s}\}\subset E \)
is a finite set containing the unit of \( E \) and such that each \( b_{i}\in B \),
\( x_{i}\in I \) and \( \Vert b_{i}\Vert ,\Vert x_{i}\Vert \leq 1 \), \( 1\leq i\leq s \)
and if \( 0\leq e\leq 1_{E} \) is an element of \( I \) such that
\( \Vert [e,\rho (b_{i})]\Vert ,\Vert x_{i}-ex_{i}\Vert ,\Vert x_{i}-x_{i}e\Vert <\delta  \)
for \( 1\leq i\leq s \) then
\[
\rcp _{E}(\omega ,30\epsilon )\leq \rcp _{I}(e^{1/n}\omega e^{1/n},\epsilon )+\rcp _{B}(\pi (\omega ),\epsilon ),\]
for any \( n>1/\epsilon  \).  
\end{lem}
\begin{proof}
Fix \( \epsilon >0 \).  By the lemma in \cite[page 332]{arveson:extensions}
we can find a \( \tilde{\delta }>0 \) such that for every pair of elements
\( f,g \) in the unit ball of \( E \) such that \( f\geq 0 \) we have the
implication
\[
\Vert [f,g]\Vert <\tilde{\delta }\Rightarrow \Vert [f^{1/2},g]\Vert <\epsilon .\]
Let \( \delta =\min (\tilde{\delta },\epsilon ^{2}/4) \). 

So assume that \( \omega =\{\rho (b_{1})+x_{1},\ldots ,\rho (b_{s})+x_{s}\}\subset E \)
is a finite set such that each \( b_{i}\in B \), \( x_{i}\in I \) and \( \Vert b_{i}\Vert ,\Vert x_{i}\Vert \leq 1 \),
\( 1\leq i\leq s \) and \( 0\leq e\leq 1_{E} \) is an element of \( I \)
such that \( \Vert [e,\rho (b_{i})]\Vert ,\Vert x_{i}-ex_{i}\Vert ,\Vert x_{i}-x_{i}e\Vert <\delta  \)
for \( 1\leq i\leq s \).
By our choice of \( \delta  \) we then have \( \Vert [e^{1/2},\rho (b_{i})]\Vert ,\Vert [(1_{E}-e)^{1/2},\rho (b_{i})]\Vert <\epsilon  \)
for all \( 1\leq i\leq s \).  

We also claim that \( \Vert [e^{1/2},x_{i}]\Vert <\epsilon  \) for all \( 1\leq i\leq s \).
 To see this we first note that since \( e\leq e^{1/2} \) we have that \( 1_{E}-e\geq 1_{E}-e^{1/2} \).
 Thus
\[
\Vert e^{1/2}x-x\Vert ^{2}=\Vert x^{*}(1_{E}-e^{1/2})^{2}x\Vert \leq \Vert x^{*}(1_{E}-e^{1/2})x\Vert \leq \Vert x^{*}(1_{E}-e)x\Vert \leq \Vert x-ex\Vert ,\]
for all \( x\in I \) with \( \Vert x\Vert \leq 1 \). Similarly, \( \Vert xe^{1/2}-x\Vert \leq \Vert x-xe\Vert  \)
and hence
\[
\Vert [e^{1/2},x_{i}]\Vert \leq \Vert e^{1/2}x_{i}-x_{i}\Vert +\Vert x_{i}e^{1/2}-x_{i}\Vert
\leq \Vert ex_{i}-x_{i}\Vert^{1/2} +\Vert x_{i}e-x_{i}\Vert^{1/2} <\epsilon .\]

Since \( 0\leq e\leq 1_{E} \) some routine functional calculus shows that \( \Vert e^{1/2+1/n}-e^{1/2}\Vert \leq 2/n \)
and \( \Vert e^{1+2/n}-e\Vert \leq 2/n \).  Combining these inequalities, the
inequalities in the previous paragraphs and a standard interpolation argument
we get
\begin{eqnarray*}
{}\Vert [e^{1/2+1/n},\rho (b_{i})]\Vert  & < & 4/n+\epsilon ,\\
{}\Vert [e^{1/2+1/n},x_{i}]\Vert  & < & 4/n+\epsilon ,\\
\Vert \rho (b_{i})e^{1+2/n}-\rho (b_{i})e\Vert  & \leq & 2/n,\\
\Vert x_{i}e^{1+2/n}-x_{i}e\Vert  & \leq & 2/n .
\end{eqnarray*}
These four inequalities will be needed at the end of the proof.  

Assume that \( E\subset B(H) \) and \( B\subset B(K) \).  By Arveson's extension
theorem we may assume that \( \rho :B\to E\subset B(H) \) is defined on all
of \( B(K) \) (and takes values in \( B(H) \)).  Now choose \( (\tilde{\phi }_{1},\tilde{\phi }_{2},D_{1})\in CPA(B) \)
such that \( \Vert \tilde{\phi }_{2}\circ \tilde{\phi }_{1}(b_{i})-b_{i}\Vert <\epsilon  \),
\( 1\leq i\leq s \), and \( \rank (D_{1})=\rcp _{B}(\pi (\omega ),\epsilon ). \)
Using the techniques in the proof of Proposition 1.4 in \cite{brown:entropy}
we may replace the (not necessarily unital) maps \( \tilde{\phi }_{1} \) and
\( \tilde{\phi }_{2} \) with unital maps \( \phi _{1}:B\to D_{1} \) and \( \phi _{2}:D_{1}\to B(K) \)
such that \( \Vert \phi _{2}\circ \phi _{1}(b_{i})-b_{i}\Vert <14\epsilon  \).
(Here we use the facts that \( 1_{B}\in \pi (\omega ) \) and \( \Vert b_{i}\Vert \leq 1,1\leq i\leq s \).) 

Similarly, choose \( (\psi _{1},\psi _{2},D_{2})\in CPA(I) \) such that
\[
\Vert \psi _{2}\circ \psi _{1}(e^{1/n}(b_{i}+x_{i})e^{1/n})-e^{1/n}(b_{i}+x_{i})e^{1/n}\Vert <\epsilon ,\]
\( 1\leq i\leq s \), and \( \rank (D_{2})=\rcp _{I}(e^{1/n}\omega e^{1/n},\epsilon ). \)
By Arveson's extension theorem we may assume that \( \psi _{1} \) is defined
on all of \( B(H) \). 

Define \( \chi _{1}:E\to D_{1}\oplus D_{2} \) by
\[
\chi _{1}(y)=\phi _{1}(\pi (y))\oplus \psi _{1}(e^{1/n}ye^{1/n}),\]
for all \( y\in E \) and \( \chi _{2}:D_{1}\oplus D_{2}\to B(H) \) by
\[
\chi _{2}(S\oplus T)=(1_{E}-e)^{1/2}\rho (\phi _{2}(S))(1_{E}-e)^{1/2}+e^{1/2}\psi _{2}(T)e^{1/2},\]
for all \( S\in D_{1},T\in D_{2} \).  Since we have arranged that \( \phi _{2} \)
is unital (and \( \rho  \) is unital by assumption) we see that \( \chi _{2}(1_{D_{1}}\oplus 1_{D_{2}})=1_{E}-e+e^{1/2}\psi _{2}(1_{D_{2}})e^{1/2}=1_{E}-\big (e^{1/2}(1_{E}-\psi _{2}(1_{D_{2}}))e^{1/2}\big ) \).
 Since \( \psi _{2} \) is a contractive completely positive map, this shows
that \( \chi _{2}(1_{D_{1}}\oplus 1_{D_{2}}) \) is a positive operator of norm
less than or equal to one.  Since it is clear that \( \chi _{2} \) is a completely
positive map this, in turn, implies that \( \chi _{2} \) is also a contractive
map (see \cite[Proposition 3.5]{paulsen:cbmaps}). 

Hence \( (\chi _{1},\chi _{2},D_{1}\oplus D_{2})\in CPA(E) \) and
\[
\rank (D_{1}\oplus D_{2})=\rcp _{I}(e^{1/n}\omega e^{1/n},\epsilon )+\rcp _{B}(\{\pi (\omega )\},\epsilon ).\]
Thus we only have to check that \( \Vert \chi _{2}\circ \chi _{1}(\rho (b_{i})+x_{i})-\rho (b_{i})-x_{i}\Vert <31\epsilon  \),
for \( 1\leq i\leq s \).  But, letting \( y=\rho (b_{i})+x_{i} \) we have 

\begin{align*}
  \|\chi_{2} \circ\chi_{1} (y) - y \|
    &=\|(1_{E} - e)^{1/2}\rho(\phi_{2}(\phi_{1}(b_{i})))(1_{E} - e)^{1/2} + \\[2mm]
    &\mathrel{\phantom\sup}\ \  e^{1/2}\psi_{2}\circ\psi_{1} 
        (e^{1/n}(\rho(b_{i}) + x_{i}) 
        e^{1/n})e^{1/2} - \rho(b_{i}) - x_{i} \|\\[2mm]
    &\leq\|(1_{E} - e)^{1/2}\rho(\phi_{2}(\phi_{1}(b_{i})))(1_{E} - e)^{1/2}\\[2mm]
    &\mathrel{\phantom\sup}\ \
        - \ (1_{E} - e)^{1/2}\rho(b_{i})(1_{E} - e)^{1/2}\|\\[2mm] 
    & + \|e^{1/2}\psi_{2}\circ\psi_{1} (e^{1/n}(\rho(b_{i}) + x_{i}) e^{1/n}) 
          e^{1/2}\\[2mm]
    &\mathrel{\phantom\sup}\ \ - \ 
        e^{1/2}\big( e^{1/n}(\rho(b_{i}) + x_{i}) e^{1/n}\big)e^{1/2}\|\\[2mm] 
    & + \| (1_{E} - e)^{1/2}\rho(b_{i})(1_{E} - e)^{1/2} + 
         e^{\frac{n+2}{2n}}(\rho(b_{i}) + x_{i}) e^{\frac{n+2}{2n}}\\[2mm]
    &\mathrel{\phantom\sup}\ \ - \rho(b_{i}) - x_{i} \| \displaybreak[2] \\[2mm] 
    &\leq15\epsilon\! + \!\| [(1_{E} - e)^{1/2}, \rho(b_{i})] \|\! + \!\| [ e^{\frac{n+2}{2n}}, \rho(b_{i})] \|\! + \!\| 
      [ e^{\frac{n+2}{2n}}, x_{i}] \|\\[2mm]
    &\mathrel{\phantom\sup}\ \ + \|\rho(b_{i})(1_{E} - e) + 
      (\rho(b_{i}) + x_{i})e^{1 + 2/n} 
      - \rho(b_{i}) - x_{i} \| \displaybreak[2] \\[2mm]
    &\leq18\epsilon+ 8/n + \|\rho(b_{i}) e^{1 + 2/n} - \rho(b_{i}) e \|\\[2mm]
    &\mathrel{\phantom\sup} + \| x_{i} e^{1 + 2/n} - x_{i} e \|\\[2mm]
    &\leq18\epsilon+ 12/n. 
\end{align*}

Hence for any \( n>1/\epsilon  \) we have the desired inequality. 
\end{proof}
\begin{rem}
The previous lemma is easily generalized to arbitrary extensions, though a precise
formulation is somewhat awkward (and does not appear to be useful for entropy
calculations).  The idea is that if \( E \) is a unital exact \( C^{*} \)--algebra
then the quotient map \( \pi :E\to B \) is always locally liftable (cf. \cite[Prop. 7.2]{kirchberg:exactness}).
 Hence the \( \delta  \)--rank of any finite subset of \( E \) can be estimated
in terms of finite subsets of \( I \) and \( B \).
\end{rem}
Though it will not be needed in what follows it seems appropriate to point out
the following application.  

\begin{prop}
With \( I,E,B,\pi  \) and \( \rho  \) as above, let \( \alpha \in \Aut (E) \)
be an automorphism such that \( \alpha (I)=I \) and let \( \hat{\alpha }\in \Aut (B) \)
be the induced automorphism.  Assume that \( \rho \circ \hat{\alpha }=\alpha \circ \rho  \)
and there exists an approximate unit \( \{e_{\lambda }\}\subset I \) such that
\( \alpha (e_{\lambda })=e_{\lambda } \), for all \( \lambda  \) (which happens,
for example, if there exists a strictly positive element \( h\in I \) such
that \( \alpha (h)=h \)).  Then
\[
\HT (\alpha )=\max (\HT (\alpha |_{I}),\HT (\hat{\alpha })).\]

\end{prop}
\begin{proof}
By \cite[Prop. 2.10]{brown:entropy} it suffices to show the inequality \( \HT (\alpha ) \leq \max (\HT (\alpha |_{I}),\HT (\hat{\alpha })) \).


Let \( \epsilon >0 \) be given and \( \omega =\{\rho (b_{1})+x_{1},\ldots ,\rho (b_{s})+x_{s}\}\subset E \)
be any finite set containing the unit of \( E \) and such that \( \Vert b_{i}\Vert ,\Vert x_{i}\Vert \leq 1 \),
\( 1\leq i\leq s \). Let \( \{e_{\lambda }\}\subset I \) be an approximate
unit such that \( \alpha (e_{\lambda })=e_{\lambda } \), for all \( \lambda  \).
 Since we can manufacture a \emph{quasicentral} approximate unit out of the
convex hull of \( \{e_{\lambda }\} \) \cite{arveson:extensions}, we may further assume that \( \{e_{\lambda }\} \)
is quasicentral in \( E \) (and still fixed by \( \alpha  \)).  

Choose \( \delta >0 \) according to the previous lemma and take \( \lambda  \)
large enough that
\( \Vert [e_{\lambda }, \rho (b_{i})]\Vert\), \(\Vert x_{i}-e_{\lambda }x_{i}\Vert,\Vert ,x_{i}-x_{i}e_{\lambda }\Vert <\delta  \)
for \( 1\leq i\leq s \).  Since \( \alpha ^{j}(e_{\lambda })=e_{\lambda } \)
for all \( j\in {\mathbb N} \) it is clear that
\( \Vert [e_{\lambda },\alpha ^{j}(\rho (b_{i}))]\Vert ,
\Vert \alpha ^{j}(x_{i})-e_{\lambda }\alpha ^{j}(x_{i})\Vert ,\Vert \alpha ^{j}(x_{i})-\alpha ^{j}(x_{i})e_{\lambda }\Vert <\delta  \)
for \( 1\leq i\leq s \) and all \( j\in {\mathbb N} \).  Hence letting \( \omega _{I}=e_{\lambda }^{1/k}\omega e_{\lambda }^{1/k}\subset I \),
for some \( k>1/\epsilon  \), the previous lemma implies that \( \rcp (\omega \cup \ldots \cup \alpha ^{n}(\omega ),30\epsilon ) \)
is bounded above by
\[
\rcp _{I}(\omega _{I}\cup \ldots \cup \alpha ^{n-1}(\omega _{I}),\epsilon )+\rcp _{B}(\pi (\omega )\cup \ldots \cup \hat{\alpha }^{n-1}(\pi (\omega )),\epsilon )\]
which is bounded above by
\[
2\max (\rcp _{I}(\omega _{I}\cup \ldots \cup \alpha ^{n-1}(\omega _{I}),\epsilon ),\rcp _{B}(\pi (\omega )\cup \ldots \cup \hat{\alpha }^{n-1}(\pi (\omega )),\epsilon )).\]
This inequality implies the result.  
\end{proof}

The next lemma is inspired by section~5 in~\cite{dvv:qcni}.

\begin{lem}
\label{lem:blktridiag}
Let $F\subset B(H)$ be a finite set
of self-adjoint contractive operators on a Hilbert space $H$.
Let $P$ be a projection in $B(H)$, of rank $\ell<\infty$. Then for any
$\delta >0$, there exists a positive finite--rank contraction $X\in B(H)$,
so that:
\begin{enumerate}
\item $XP=PX=P$ 
\item $\Vert [X,T]\Vert <\delta $ for all $T\in F$ 
\item The rank of $X$ is no bigger than $\ell\cdot(|F|+1)^{(2/\delta )+1}$.
\end{enumerate}
\end{lem}
\begin{proof}
Denote by $K_{1}\subset H$ the range of $P$. Define recursively
\[
K_{n}=\textrm{span}\{K_{n-1}\cup \bigcup _{T\in F}TK_{n-1}\}.
\]
Let $q=|F|$ be the cardinality of $F$.
Then $K_{n}$ has dimension at
most $(q+1)$ times the dimension of $K_{n-1}$, so that $\dim K_{n}\leq\ell\cdot(q+1)^{n}$. 

Let $P_{n}$ be the orthogonal projection onto $K_{n}$; then $P_{n}$
are clearly an increasing sequence, and $P_{i}P=PP_{i}=P=P_{1}$ for all
$i$. Let
\[
X_{n}=\frac{1}{n}(P_{1}+\cdots +P_{n}).\]
Then $X_{n}P=PX_{n}=\frac{1}{n}(PP_{1}+\cdots +PP_{n})=P$ for all $n$.
Note that the rank of $X_{n}$ is the same as that of $P_{n}$, which
is bounded by $\ell\cdot(q+1)^{n}$ 

Set $Q_{n}=P_{n}-P_{n-1}$, $Q_{1}=P_{1}$.
Since $TK_{n}\subset K_{n+1}$ and thus
$TP_{n}=P_{n+1}TP_{n}$,
if $m-n>2$ then we have $Q_{m}TQ_{n}=Q_{m}TP_{n}Q_{n}=Q_{m}P_{n+1}TP_{n}Q_{n}=0$ for all $T\in F$.
Since $T$ is self-adjoint,
also $Q_{n}TQ_{m}=0$ if $m-n>2$. Hence $Q_{n}TQ_{m}=0$ if $|n-m|>2$.
Let $Z$ be the orthocomplement of $\overline{\cup K_{n}}$ in $H$.
Then in the decomposition $H=\bigoplus Q_{n}H\oplus Z$,
each $T\in F$ has the form
\[
\left( \begin{array}{ccccc}
Q_{1}TQ_{1} & Q_{1}TQ_{2} & 0 &  & \\
Q_{2}TQ_{1} & Q_{2}TQ_{2} & Q_{2}TQ_{3} & 0 & \\
0 & Q_{3}TQ_{2} & Q_{3}TQ_{3} & Q_{3}TQ_{4} & \ddots \\
 & 0 & \ddots  & \ddots  & \ddots \\
 &  & \ddots  & \ddots  & \quad 0\quad 
\end{array}\right) ,\]
i.e., it is a ``block tri--diagonal'' matrix.
Hence using the convention $Q_0=0$ we have the identities
\begin{eqnarray*}
Q_{n}T=\sum _{j=-1,0,1}Q_{n}TQ_{n+j}, &  & \\
TQ_{n}=\sum _{j=-1,0,1}Q_{n+j}TQ_{n}. &  & 
\end{eqnarray*}

Now $X_{n}=Q_{1}+(1-\frac{1}{n})Q_{2}+\cdots +\frac{1}{n}Q_{n}=\sum _{i=1}^{n}(\frac{n-i+1}{n})Q_{i}$
so that for any $T\in F$,
\[
X_{n}T=\sum _{i=1}^{n}\biggl(\frac{n-i+1}{n}\biggr)Q_{i} T=\sum_{j=-1,0,1}\sum_{i=1}^{n}\biggl(\frac{n-i+1}{n}\biggr)Q_{i}TQ_{i+j}.
\]
Similarly,
\begin{eqnarray*}
TX_{n}=\sum _{i=1}^{n}\biggl(\frac{n-i+1}{n}\biggr)TQ_{i}& = & \sum _{j=-1,0,1}\sum _{i=1}^{n}\biggl(\frac{n-i+1}{n}\biggr)Q_{i+j}TQ_{i}\\
 & = & \sum _{j=-1,0,1}\sum _{i=1}^{n}\biggl(\frac{n-i+1-j}{n}\biggr)Q_{i}TQ_{i+j}
\end{eqnarray*}
Hence
\begin{eqnarray*}
\Vert X_{n}T-T_{n}X\Vert  & = & \bigg\Vert \sum _{j=-1,0,1}\sum _{i=1}^{n}\biggl(\frac{(n-i+1)-(n-i+1-j)}{n}\biggr)Q_{i}TQ_{i+j}\bigg\Vert \\
 & = & \bigg\Vert \sum _{j=-1,0,1}\sum _{i=1}^{n}\frac{j}{n}Q_{i}TQ_{i+j}\bigg\Vert \\
 & \leq  & \sum _{j=-1,0,1}\frac{|j|}{n}\bigg\Vert \sum _{i}Q_{i}TQ_{i+j}\bigg\Vert \\
 & \leq  & \sum _{j=-1,1}\frac{|j|}{n}=\frac2n
\end{eqnarray*}
The last inequality is due to the fact that for a fixed $j$, the operator
$\sum Q_{i}TQ_{i+j}$ is a block-diagonal operator, with the blocks having
orthogonal ranges, so that $\Vert \sum Q_{i}TQ_{i+j}\Vert \leq \max \{\Vert Q_{i}TQ_{i+j}\Vert \}\leq 1$. 

Choose the smallest integer $n$ with $n>\frac{2}{\delta }$, and set
$X=X_{n}$. Then $\Vert [X,T]\Vert \leq \frac{2}{n}<\delta $; by construction,
$XP=PX=P$ and the rank of $X$ is bounded by $\ell(q+1)^{n}\leq\ell(q+1)^{\frac{2}{\delta }+1}$.
\end{proof}

We finally come to the main result of this section.  

\begin{thm}
\label{theorem:entropyextensionestimate} Let $A,B\subset B(H)$ be unital
exact $C^{*}$--algebras such that $B\cap \mathcal{K}(H)=\{0\}$ and
$C=B\otimes1+\mathcal{K}(H)\otimes A\subset B(H\otimes H)$.  For any two
unitaries $V,W\in B(H)$ such that ${\Ad }V(B)=B$ and ${\Ad }W(A)=A$,
the unitary $U=V\otimes W$ has the property that ${\Ad }U(C)=C$ and
\[
\HT ({\Ad }U|_C)\leq \max (\HT ({\Ad }V|_B),\HT ({\Ad }W|_A)).\]
\end{thm}
\begin{proof}
We have a short exact sequence
\begin{equation}
\label{eq:AKCB}
0\to\mathcal{K}(H)\otimes A\to C\overset\pi\to B\to0
\end{equation}
with splitting $\rho:B\to C$ given by $\rho(b)=b\otimes1$.

Let $k_{1},\dots ,k_{N}\in \mathcal{K}(H)$, $a_{1},\dots ,a_{N}\in A$ and $b_{1},\dots ,b_{N}\in B$
be self-adjoint elements, each of norm at most $1$, so that each $k_{i}$ has finite rank.
Let $L$ be the sum of the ranks of $k_{1},\dots ,k_{N}$.
Let
\[
\omega=\{k_i\otimes a_i+b_i\otimes1:1\le i\le N\}.
\]
Fixing $\epsilon>0$ and a positive integer $n$, let
\begin{eqnarray*}
\omega(n) & = & \omega\cup(\Ad U)(\omega)\cup\cdots\cup(\Ad U)^{n-1}(\omega) \\
          & = & \{1\}\cup\{V^rk_iV^{-r}\otimes W^ra_iW^{-r}+V^rb_iV^{-r}\otimes1:1\le i\le N,\,0\le r\le n-1\}.
\end{eqnarray*}
Then the
sum of the ranks of
\[
k_{1},\dots ,k_{N},\dots ,V^{n-1}k_{1}V^{-(n-1)},\dots ,V^{n-1}k_{N}V^{-(n-1)}\]
is at most $nL$, so there exists a projection $P$ of rank $nL$,
so that $PV^{r}k_{j}V^{-r}=V^{r}k_{j}V^{-r}P=V^{r}k_{j}V^{-r}$ for all
$j=1,\dots ,N$ and $r=0,\dots ,n-1$.

Consider the collection
\[
F=\{b_{1},\dots ,b_{N}\}\cup \cdots \cup \{V^{n-1}b_{1}V^{-(n-1)},\dots ,V^{n-1}b_{N}V^{-(n-1)}\},
\]
of at most $nN$ self-adjoint operators on $B(H)$.
Let $\delta $ be as in Lemma \ref{label:rcpestimateextension} for the given value of $\epsilon$ and the short exact sequence~\eqref{eq:AKCB}.
By applying Lemma \ref{lem:blktridiag} with these choices of $\delta $,
$P$, $F$, and $\ell=nL$,
we find a positive, finite--rank $X\in \mathcal{K}(H)$, so that:
\begin{enumerate}
\item $XP=PX=P$, hence $XV^{r}k_{j}V^{-r}=XPV^{r}k_{j}V^{-r}=PV^{r}k_{j}V^{-r}=V^{r}k_{j}V^{-r}=V^{r}k_{j}V^{-r}X$
for all $1\le j\le N$ and $0\leq r\leq n-1$
\item $\Vert [X,V^{r}b_{j}V^{-r}]\Vert \leq \delta $ for all $1\le j\le N$ and
$0\leq r\leq n-1$
\item The rank of $X$ is at most $nL\cdot (nN+1)^{(2/\delta )+1}$.
\end{enumerate}
Let $e=X\otimes1$.
Then $e$ satisfies the hypotheses of Lemma
\ref{label:rcpestimateextension}.
Indeed,
\[
\Vert [e,V^{r}b_{i}V^{-r}\otimes1]\Vert =\Vert [X,V^{r}b_{i}V^{-r}]\Vert \leq \delta ,\quad 1\leq i\leq N,\,0\leq r\leq n-1
\]
and
\begin{align*}
eV^{r}k_{i}V^{-r}\otimes W^{r}a_{i}W^{-r}-V^{r}k_{i}V^{-r}\otimes W^{r}a_{i}W^{-r}&= \\
=(XV^{r}k_{i}V^{-r}-V^{r}k_{i}V^{-r})\otimes W^{r}a_{i}W^{-r}&=0
\end{align*}
if $1\leq i\leq N$, and $0\leq r\leq n-1$, etc. Hence by Lemma~\ref{label:rcpestimateextension},
as long as $m>\epsilon^{-1}$ we have
\[
\rcp_{C}(\omega(n),30\epsilon)\le\rcp_{\mathcal{K}(H)\otimes A}(e^{1/m}\omega (n)e^{1/m},\epsilon )+\rcp _{B}(\pi (\omega (n)),\epsilon ).
\]
Note that
\begin{align*}
e^{1/m}\omega (n)e^{1/m}=&\{X^{2/m}\otimes 1\}\cup \\
&\cup
\begin{aligned}[t]
\{X^{1/m}V^rk_iV^{-r}X^{1/m}\otimes W^ra_iW^{-r}+&X^{1/m}V^rb_iV^{-r}X^{1/m}\otimes1 \\
&:1\le i\le N,\,0\le r\le n-1\}.
\end{aligned}
\end{align*}
Hence setting
\begin{align*}
\omega _{A}(n)&=\{1\}\cup\{W^ra_{i}W^{-r}:1\le i\le N,\,0\le r\le n-1\} \\
\omega _{B}(n)&=\{1\}\cup\{V^rb_{i}V^{-r}:1\le i\le N,\,0\le r\le n-1\}
\end{align*}
we have
\[
\rcp _{\mathcal{K}(H)\otimes A}(e^{1/n}\omega (n)e^{1/n},\epsilon )\leq \rcp _{A}(\omega _{A}(n),\epsilon/2 )\cdot \textrm{rank}(X).
\]
On the other hand,
\[
\rcp _{B}(\pi (\omega (n)),\epsilon )=\rcp _{B}(\omega _{B}(n),\epsilon ).\]
Since the rank of $X$ is at most $nL(nN+1)^{(\delta /2)+1}$
it follows that
\begin{eqnarray*}
\rcp(\omega(n),30\epsilon)&\le&\rcp_{A}(\omega_{A}(n),\epsilon/2)nL(nN+1)^{2/\delta +1}
+\rcp_{B}(\omega_{B}(n),\epsilon) \\
 & \leq & 2(nL+1)(nN+1)^{2/\delta +1}\max\big(\rcp_{A}(\omega _{A}(n),\epsilon/2),\rcp_{B}(\omega_{B}(n),\epsilon)\big)
\end{eqnarray*}
so that
\begin{align*}
\limsup_{n\to\infty}\frac{1}{n}&\log\rcp(\omega(n),30\epsilon)\quad\le \\
&\le\quad \max\bigg(\limsup_{n\to\infty}\frac{1}{n}\log\rcp(\omega_{A}(n),\epsilon/2),\,\limsup_{n\to\infty}\frac{1}{n}\log\rcp(\omega_{B}(n),\epsilon)\bigg)\\
 & \le\quad\max\big(\HT (\Ad W|_{A}),\HT (\Ad V|_{B})\big)_.
\end{align*}
Since $C$
can be written as the closure of the linear span of elements of the form appearing in $\omega$,
the statement of the theorem follows. 
\end{proof}

We shall record the following corollary, which will be the basis for entropy
computations in this paper.

\begin{cor}
\label{Corr:algebrasAn}Let \( A \) be a \( C^{*} \)--algebra, \( \alpha \in \Aut (A) \),
\( \pi :A\to B(H) \) be a faithful representation of \( A \) and \( U\in U(H) \)
be such that \( \pi (\alpha (a))=U\pi (a)U^{*} \). Assume that \( \pi (A)\cap K(H)=\{0\} \).
Let \( A_{0}=A \), \( \alpha _{0}=\alpha  \), \( \pi _{0}=\pi  \), \( H_{0}=H \)
and \( U_{0}=U \).
Recursively construct $C^*$--algebras \( A_{n},\alpha _{n}\in \Aut (A_{n}) \),
\( \pi _{n}:A_{n}\to B(H_{n}) \) and \( U_{n}\in U(H_{n}) \) by setting
\begin{eqnarray*}
H_{n} & = & H_{n-1}\otimes H,\\
A_{n} & = & K(H_{n-1})\otimes \pi (A)+\pi _{n-1}(A_{n-1})\otimes \Id _{H},\\
\pi _{n} & = & \text {obvious representation on }H_{n},\\
U_{n} & = & U_{n-1}\otimes U,\\
\alpha _{n} & = & \Ad U_{n}.
\end{eqnarray*}
View \( A_{n-1}\subset A_{n} \) as \( A_{n-1}\cong \pi _{n-1}(A_{n})\otimes \Id _{H} \).
Let
\[
A_{\infty }=\overline{\cup _{n}A_{n},}\qquad \alpha _{\infty }=\lim _{\to }\alpha _{n}.\]
Then:

(i) \( \HT (\alpha _{\infty })=\HT (\alpha ) \). 

(ii) If \( \gamma  \) is an injective endomorphism of \( A_{\infty } \), so that \( \gamma \circ \alpha _{\infty }=\alpha _{\infty }\circ \gamma  \),
denote by \( \bar{\alpha } \) the obvious extension of \( \alpha _{\infty } \)
to \( A_{\infty }\rtimes _{\gamma }\mathbb {N} \). Then \( \HT (\bar{\alpha })=\HT (\alpha ) \).
\end{cor}
\begin{proof}
Statement (ii) follows from statement (i) and the results of \cite{shlyakht-dykema:exactness}.
Hence it is sufficient to prove (i); for that one only needs to prove that \( \HT (\alpha _{n})\leq \HT (\alpha ) \),
in view of the behavior of entropy with respect to inductive limits. We now proceed
by induction on \( n \). Since \( \alpha _{0}=\alpha  \), the statement is
true for \( n=0 \). Applying Theorem \ref{theorem:entropyextensionestimate}
to \( A_{n}=K(H_{n-1})\otimes \pi (A)+\pi _{n-1}(A_{n-1})\otimes \Id _{H} \)
gives \( \HT (\alpha _{n})\leq \max (\HT (\alpha _{n-1}),\HT (\alpha )) \),
which is equal to \( \HT (\alpha ) \) by the induction hypothesis.
\end{proof}

\section{Free products with the Toeplitz algebra.}

The main technical result of this section (which will be used to prove the more
general result about entropy of amalgamated free products of automorphisms)
is the following theorem.

\begin{thm}
\label{Theorem:freeprodToeplitz}Let \( \alpha  \) be an automorphism of a
\( D \)--probability space \( (A,E:D\to A) \). Assume that \( D \) is finite dimensional.
Assume that the GNS representation associated to \( E \) is faithful. Let \( \mathcal{T}\subset B(\ell ^{2}) \)
be the Toeplitz algebra generated by the unilateral shift \( \ell (\delta _{n})=\delta _{n+1} \)
(\( n\geq 1 \)), and \( \psi  \) be the vector state on \( \mathcal{T} \)
associated to \( \delta _{1}\in \ell ^{2} \). Consider on the algebra \( (A,E)*_{D}(D\otimes \mathcal{T},\id _{D}\otimes \psi ) \)
the automorphism \( \alpha *(\alpha |_{D}\otimes \id ) \). Then \( \HT (\alpha *(\alpha |_{D}\otimes \id )) = \HT (\alpha )\).
\end{thm}
Because of Corollary \ref{Corollary:embeddingToeplitz} and monotonicity of
\( \HT  \), it is sufficient to prove Theorem \ref{Theorem:freeprodToeplitz}
in the particular case that \( D=\mathbb {C} \). For convenience, we'll restate
this particular case as 

\begin{prop}
Let \( A \) be a unital \( C^{*} \)--algebra, \( \phi :A\to \mathbb {C} \)
a state with a faithful GNS representation. Let \( \alpha \in \Aut (A) \) be
an automorphism, so that \( \phi \circ \alpha =\phi  \) . Consider the algebra
\( (A,\phi )*(\mathcal{T},\psi ) \). Then \( \HT (\alpha *\id )=\HT (\alpha ) \).
\end{prop}
\begin{proof}
Let \( H=L^{2}(A,\phi ) \) be the GNS Hilbert space associated to \( A
\), and let \( 1\in H \) be the cyclic vector associated to \( \phi
\). Let \( \mathcal{O}_{2} \) be the Cuntz algebra on two generators
\cite{cuntz}. Without loss of generality, by replacing \( A \) with \(
A\otimes \mathcal{O}_{2} \) and \( \alpha \) with \( \alpha \otimes
\id \) we may assume that the GNS representation \( \pi :A\to B(H) \)
satisfies \( \pi (A)\cap K(H)=\{0\} \). Let \( U:H\to H \) be the
unitary induced on \( H=L^{2}(A,\phi ) \) by \( \alpha \). We shall
covariantly identify \( ((A,\phi )*(\mathcal{T},\psi ),\alpha *\id )
\) as the crossed product by an certain endomorphism \( \gamma \) of
the algebra \( A_{\infty } \) described in Corollary
\ref{Corr:algebrasAn}, taken with the automorphism \( \bar{\alpha }
\).  By Corollary \ref{Corr:algebrasAn}, we then have \( \HT (\alpha
)=\HT (\bar{\alpha }) \), which, in view of our identification, is the
same as \( \HT (\alpha *\id ) \), hence proving the Proposition.  The
remainder of the proof is essentially a special case of the techniques
used in \cite{shlyakht-dykema:exactness} where more general Cuntz--Pimsner algebras were
shown to have a crossed product structure.

Consider now the Hilbert space
\[
F=L^{2}(A,\phi )\oplus \bigoplus _{n\geq 2}L^{2}(A,\phi )^{\otimes n},\]
and the representation \( \rho :A\to B(H) \) given by \( \rho =\pi \oplus (\bigoplus _{n\geq 2}\pi \otimes 1) \).
Consider the isometry \( \ell :F\to F \) defined by
\[
\ell (\xi )=1\otimes \xi ,\]
where \( 1 \) denotes the image of the unit of \( A \) in \( L^{2}(A,\phi ) \),
and \( \xi \in F \). Denote by \( U:L^{2}(A,\phi )\to L^{2}(A,\phi ) \) the
unitary implementing \( \alpha  \). Denote by \( V \) the unitary \( U\oplus \bigoplus _{n}U^{\otimes n} \)
acting on the Hilbert space \( F \). One sees that \( \Ad _{V}(\ell )=\ell  \),
and \( \Ad _{V}(\rho (a))=\rho (\alpha (a)) \). 
\begin{claim}
\( (C^{*}(\rho (A),\ell ),\Ad _{V})\cong ((A,\phi )*(\mathcal{T},\psi ),\alpha *\id ) \). 
\end{claim}
\begin{proof}
This actually follows from \S\ref{sec:exampleBogoljubov}, since \( C^{*}(\rho (A),\ell ) \)
is isomorphic to the Cuntz--Pimsner algebra associated to the \( A,A \) bimodule
\( L^{2}(A,\phi )\otimes A \). For the reader's convenience, we give a proof.

Let \( \theta  \) be the vector--state on \( B(F) \), associated to the vector
\( 1\in L^{2}(A,\phi )\oplus 0\subset F \). Then \( \ell ^{*}1=0 \), and one
can easily verify that: (i) \( \theta (a_{0}\ell \cdots \ell a_{n}\ell ^{*}a_{n+1}\cdots \ell ^{*}a_{n+m})=0 \)
for all \( a_{j}\in \rho (A) \), \( n,m\geq 0 \), \( n+m>0 \); and (ii) \( \ell ^{*}\rho (a)\ell =\phi (a) \),
for all \( a\in A \). It follows from \cite{shlyakht:amalg} that \( A \)
and \( C^{*}(\ell ) \) are free in \( (B(F),\theta ) \). Since \( C^{*}(\ell ) \)
is (obviously) isomorphic to \( \mathcal{T} \), \( \theta |_{C^{*}(\ell )}=\psi  \),
and since \( \rho  \) is injective and \( \theta |_{\rho (A)}=\phi  \), the
claim is proved. 
\end{proof}
 From now on, write \( C=C^{*}(A,\ell ) \). 

Denote by \( C_{n} \) the closed linear span
\[
C_{n}=\overline{\Span} \{a_{0}\ell a_{1}\cdots \ell a_{m}\ell ^{*}a_{m+1}\ell ^{*}
\cdots a_{2m}:m\leq n,\, \, a_{0},\dots ,a_{2m}\in \rho (A)\}.\]
(each monomial above has exactly \( m \) terms equal to \( \ell  \) and the
same number of terms equal to \( \ell ^{*} \)). Note that because of the relation
\( \ell ^{*}\rho (a)\ell =\phi (a) \), \( a\in A \), each \( C_{n} \) is
a \( C^{*} \)--subalgebra of \( C \). The action of \( w=a_{0}\ell a_{1}\cdots \ell a_{n}\ell ^{*}a_{n+1}\ell ^{*}\cdots a_{2n}\in C_{k} \)
on a vector \( \xi =\xi _{1}\otimes \cdots \otimes \xi _{r}\in H \) can be
described as follows:
\begin{eqnarray*}
w\cdot \xi  & = & 0,\qquad \textrm{if }n\geq r\\
w\cdot \xi  & = & (a_{0}\otimes a_{1}\otimes \cdots \otimes (a_{n}\cdot \xi _{n+1})\otimes \dots \otimes \xi _{r})\prod _{j=1}^{n}\langle \xi _{j},a_{2n-j+1}^{*}\rangle ,\, \textrm{otherwise}
\end{eqnarray*}
(here we identify \( a_{j}\in \rho (A) \) with \( a_{j}\cdot 1\in L^{2}(A,\phi )\subset H \)).

Denote by \( C_{\infty } \) the $C^*$--algebra
\( \overline{\cup _{n\geq 0}C_{n}} \). Then \( \ell C_{n}\ell ^{*}\subset C_{n+1} \)
and hence \( \gamma =\ell \cdot \ell ^{*} \) determines an endomorphism of
\( C_{\infty } \). It can be easily seen that \( C \) is isomorphic to the
crossed product of \( C_{\infty } \) by this endomorphism (cf.\ the discussion
after Proposition 1.2 in \cite{shlyakht-dykema:exactness} for the definition).
Moreover, \( \Ad _{V} \) leaves \( C_{\infty } \) invariant and commutes with
the endomorphism \( \ell \cdot \ell ^{*} \) (\cite[Proposition 1.5]{shlyakht-dykema:exactness}). 

It remains to show that \( (C_{n},\Ad _{V}|_{C_{n}})\cong (A_{n},\alpha _{n}) \),
where \( A_{n} \) and \( \alpha _{n} \) are as in Corollary \ref{Corr:algebrasAn}.

We proceed by induction. In the case that \( n=0 \), \( C_{0}=A \), \( \Ad _{V}|_{C_{0}}=\alpha  \).

Let
\[
F_{n}=L^{2}(A,\phi )\oplus \bigoplus _{k=2}^{n+1}L^{2}(A,\phi )^{\otimes k}\subset F.\]
Then \( F_{n} \) is invariant under the action of \( C_{n} \) on \( F \).
Consider the isomorphism \( W \)
\begin{eqnarray*}
W:L^{2}(A,\phi )\oplus \bigoplus _{k\geq 2}L^{2}(A,\phi )^{\otimes k}\to  &  & \\
\left( L^{2}(A,\phi )\oplus \bigoplus _{k=2}^{n+1}L^{2}(A,\phi )^{\otimes k}\right) \otimes \left( \mathbb {C}\oplus \bigoplus _{m\geq 1}L^{2}(A,\phi )^{\otimes m(n+1)}\right) =F_{n}\otimes K_{n} & 
\end{eqnarray*}
For each \( a\in C_{n} \), \( W^{*}aW:F_{n}\otimes K_{n}\to F_{n}\otimes K_{n} \)
has the form \( a|_{F_{n}}\otimes \Id _{K_{n}} \). It follows that the representation
of \( C_{n} \) on \( F_{n} \) is faithful.

Denote by \( B_{n}\subset C_n \) the ideal
\[
B_{n}=\Span \{a_{0}\ell a_{1}\cdots \ell a_{n}\ell ^{*}a_{n+1}\ell ^{*}\cdots a_{2n}:\, \, a_{0},\dots ,a_{2n}\in \rho (A)\}.\]
Write
\[
F_{n}=L^{2}(A,\phi )\oplus \bigoplus _{k=2}^{n+1}L^{2}(A,\phi )^{\otimes k}=F_{n-1}\oplus L^{2}(A,\phi )^{\otimes (n+1)}.\]
Then \( B_{n}\cdot F_{n-1}=0 \); hence the representation of \( B_{n} \) obtained
by restricting its action to the space \( L^{2}(A,\phi )^{\otimes (n+1)} \)
is faithful. The action of \( a_{0}\ell a_{1}\cdots \ell a_{n}\ell ^{*}a_{n+1}\ell ^{*}\cdots a_{2n}\in B_{n} \)
on \( \xi _{1}\otimes \cdots \otimes \xi _{n+1}\in L^{2}(A,\phi )^{\otimes (n+1)} \)
can be explicitly written as
\begin{eqnarray*}
a_{0}\ell a_{1}\cdots \ell a_{n}\ell ^{*}a_{n+1}\ell ^{*}\cdots a_{2n}\cdot \xi _{1}\otimes \cdots \otimes \xi _{n}\otimes \xi _{n+1}= &  & \\
\langle \xi _{1},a_{2n}^{*}\rangle \cdots \langle \xi _{n},a_{n+1}^{*}\rangle a_{0}\otimes \dots \otimes a_{n-1}\otimes (a_{n+1}\cdot \xi _{n+1}). &  & 
\end{eqnarray*}
Denote by \( \theta (a_{0},\dots ,a_{n-1},a_{n+1},\dots ,a_{2n})\in B(L^{2}(A,\phi )^{\otimes n}) \)
the compact operator given by
\[
\theta (a_{0},\dots ,a_{n-1}^{*},a_{n+1},\dots ,a_{2n})\xi _{1}\otimes \cdots \otimes \xi _{n}=\langle \xi _{1},a_{2n}^{*}\rangle \cdots \langle \xi _{n},a_{n+1}^{*}\rangle a_{0}\otimes \dots \otimes a_{n-1}.\]
Then the map
\begin{eqnarray*}
a_{0}\ell a_{1}\cdots \ell a_{n}\ell ^{*}a_{n+1}\ell ^{*}\cdots a_{2n}\mapsto \theta (a_{0},\dots ,a_{n-1}^{*},a_{n+1},\dots ,a_{2n})\otimes a_{n} &  & \\
\in K(L^{2}(A,\phi )^{\otimes n})\otimes A\subset B(L^{2}(A,\phi )^{\otimes n}\otimes L^{2}(A,\phi )) &  & 
\end{eqnarray*}
is a \( C^{*} \)--algebra isomorphism of \( B_{n} \) with \( K(L^{2}(A,\phi )^{\otimes n})\otimes A \).

For all \( m\leq n+1 \), the subspaces \( L^{2}(A,\phi )^{\otimes m}\subset F_{n} \)
are invariant under the action of \( C_{n} \). Denote by \( \rho _{n} \) the
representation of \( C_{n} \) obtained by restricting its action on \( H \)
to the space \( L^{2}(A,\phi )^{\otimes (n+1)} \). The image \( \rho _{n}(A_{n-1}) \)
lies inside \( B(L^{2}(A,\phi )^{\otimes (n-1)})\otimes A\otimes \Id _{L^{2}(A,\phi )} \).
By our assumption, \( \pi (A)\cap K(L^{2}(A,\phi ))=\{0\} \). This implies
that \( \rho _{n}(C_{n})\cap \rho _{n}(B_{n})=\rho _{n}(C_{n})\cap K(L^{2}(A,\phi )^{\otimes n})\otimes A=\{0\} \). 

We claim that \( \rho _{n} \) is faithful. Indeed, assume it's not, and some
\( 0\neq a\in C_{n} \) is annihilated by \( \rho  \). Write \( a=a_{n-1}+b \),
where \( a_{n-1}\in C_{n-1} \) and \( b\in B \). Since \( \rho _{n}(a)=\rho _{n}(a_{n-1})+\rho _{n}(b)=0 \),
and \( \rho _{n}(B_{n})\cap \rho _{n}(C_{n-1})=\{0\} \), it follows that \( \rho _{n}(a_{n-1})=0 \)
and \( \rho _{n}(b)=0 \). But for \( a_{n-1}\in C_{n-1} \), \( \rho (a_{n-1})=\rho _{n-1}(a_{n-1})\otimes \Id _{L^{2}(A,\phi )} \).
Proceeding inductively, we see that \( \rho _{n} \) is faithful.

We have therefore proved that
\[
C_{n}=\rho _{n-1}(A_{n-1})\otimes \Id _{L^{2}(A,\phi )}+K(L^{2}(A,\phi )^{n-1})\otimes \pi (A),\]
which means that \( C_{n}\cong A_{n} \). Since \( \Ad V_{n}|_{C_{n}} \) is
given by \( \Ad _{U^{\otimes n-1}\otimes U} \), it follows that the
dynamical system \( (C_{n},\Ad V|_{C_{n}})\cong (A_{n},\alpha _{n}) \).
\end{proof}
\begin{cor}
Let \( \alpha  \) be an automorphism of a non--commutative probability space
\( (A,E:A\to D) \), and assume that \( D \) is finite dimensional. Assume
that the GNS representation with respect to \( E \) is faithful. Let \( \bar{\alpha } \)
be the Bogoljubov automorphism of the Cuntz--Pimsner \( C^{*} \)--algebra described
in \S\ref{sec:exampleBogoljubov}. Then \( \HT (\bar{\alpha })=\HT (\alpha ) \).
\end{cor}
It would be interesting to find a formula for the entropy of more general Bogoljubov
automorphisms of Cuntz--Pimsner algebras. For example, we believe that the Corollary
above should hold for more general \( D \).

\section{Entropy for free products of automorphisms.}

We have now almost arrived at the main result of the paper 
(Theorem~\ref{Theorem:mainresult}), calculating entropy of the free product
of two automorphisms.  First, we shall prove several technical lemmas which
reduce the general case to the case of a free product with the Toeplitz algebra as in
Theorem~\ref{Theorem:freeprodToeplitz}.

\begin{lem}
\label{Lemma:freenessconjugation}
Let \( D\subset B\subset \mathcal{B} \) \( C^{*} \)--subalgebras, and let \( U\in\mathcal{B} \)
be a unitary. Let \( E \) be a conditional expectation from \( \mathcal{B} \)
onto \( D \). Assume that: (i) \( C^{*}(D,U) \)  and $B$ are free with respect to $E$
and that (ii) \( [U,D]=0 \) and \( E(U)=E(U^{*})=0 \).
Then the algebras \( B \) and \( UBU^{*} \) are free with respect to $E$.
If in addition, \( E(U^{k})=0 \) for all \( k\neq 0 \), then the
algebras \( \{U^{k}BU^{-k}\}_{k\in \mathbb {Z}} \) are free with respect to $E$.
\end{lem}
\begin{proof}
Since \( U \) is free from \( B \) and commutes with \( D \), we see that
for all \( b\in B \), \( E(UbU^{*})=E(U(b-E(b))U^{*})+E(UE(b)U^{*})=0+E(E(b)UU^{*})=E(b) \).
Now let \( b_{j}\in B \) be such that \( E(b_{j})=0 \). Then
\[
E(b_{0}Ub_{1}U^{*}b_{2}Ub_{3}U^{*}\cdots )=0\]
because \( U \) and \( B \) are free with respect to $E$. But
then \( B \) and \( UBU^{*} \) are free with respect to $E$, since
for all \( j \), \( Ub_{2j+1}U^{*}\in UBU^{*}\cap \ker E \), and any element
in \( UBU^{*}\cap \ker E \) has the form \( UbU^{*} \), for some \( b\in B\cap \ker E \).

The proof that \( \{U^{k}BU^{-k}\}_{k} \) are free with respect to $E$
proceeds along similar lines.
\end{proof}

The following lemma is included for completeness.
We will use the implications (i) $\Rightarrow$ (ii) $\Rightarrow$ (iii) in the sequel.

\begin{lem}
\label{lem:condexp}
Let $D\subseteq B\subseteq A$ be unital inclusions of unital $C^*$--algebras, and suppose there
are conditional expectations $E^A_B:A\to B$ and $E_D^B:B\to D$.
Let $E^A_D=E^B_D\circ E^A_B:A\to D$.
Consider the following statements:
\renewcommand{\labelenumi}{(\roman{enumi})}
\begin{enumerate}
\item The GNS representations associated to $E^A_B$ and $E^B_D$ are faithful.
\item The GNS representation associated to $E^A_D$ is faithful.
\item The GNS representation associated to $E^A_B$ is faithful.
\end{enumerate}
Then we have (i)$\Rightarrow$(ii)$\Rightarrow$(iii).
None of the reverse implications hold in general.
\end{lem}
\begin{proof}
We will denote by $a\to\hat a$ the defining map from $A$ onto a dense subspace of $L^2(A,E^A_D)$, and
similarly for the other $L^2$ spaces.
We will also need the isomorphism $\pi:L^2(A,E^A_D)\to L^2(A,E^A_B)\otimes_BL^2(B,E^B_D)$ given by
$\pi(\hat a)=\hat a\otimes\hat1$.

(i)$\Rightarrow$(ii): Let $a\in A$.  Choose $a_1\in A$ so that
$(aa_1)\hat{\;}\in L^2(A,E^A_B)\backslash\{0\}$,
i.e.\ $E^A_B(a_1^*a^*aa_1)$ $\ne0$.
There exists and element  $b_1\in B$ such that $$(E^A_B(a_1^*a^*aa_1))^{1/2}b_1\in L^2(B,E^B_D)\backslash\{0\},$$
i.e.\ $$0\ne E^B_D(b_1^*E^A_B(a_1^*a^*aa_1)b_1)=E^A_D(b_1^*a_1^*a^*aa_1b_1).$$
Thus $(aa_1b_1)\hat{\;}\in L^2(A,E^A_D)\backslash\{0\}$.

(ii)$\Rightarrow$(iii): Given $a\in A$ take $a_1\in A$ such that $E^A_D(a_1^*a^*aa_1)\ne0$, hence
$E^A_B(a_1^*a^*aa_1)\ne0$.

A counter example to (ii)$\Rightarrow$(i) is provided by taking
$D=\mathbb{C}$, $B=\mathbb{C}\oplus\mathbb{C}$, $A=M_2(\mathbb{C})$,
$E^A_B((c_{ij})_{1\le i,j\le2})=c_{11}\oplus c_{22}$ and $E^B_D(z_1\oplus z_2)=z_1$.

A counter example to (iii)$\Rightarrow$(ii) is provided by taking $D$, $B$ and $E^B_D$
as in the above example,
$A=M_2(\mathbb{C})\oplus M_2(\mathbb{C})$ and
$E^A_B((c_{ij})_{1\le i,j\le2}\oplus(d_{ij})_{1\le i,j\le2})=c_{11}\oplus d_{11}$.
\end{proof}



\begin{lem} \label{lemma:embedfreeprodintodirectsum}
Let \( E_{i}:A_{i}\to D \), \( i=1,2 \) be \( D \)--probability spaces, with
automorphisms \( \alpha _{i} \). Assume that the GNS representations associated
to \( E_{i} \) are faithful, and that \( \alpha _{1}|_{D}=\alpha _{2}|_{D} \).
Let \( E_{1}\oplus E_{2}:A_{1}\oplus A_{2}\to D\oplus D \) be the obvious conditional
expectation. Consider the algebra \( M_{2}(D) \) of \( 2\times 2 \) matrices
over \( D \), and view \( D\oplus D\subset M_{2}(D) \) as diagonal matrices.
Let \( F:M_{2}(D)\to D\oplus D \) be the conditional expectation \( F\left( \begin{smallmatrix}
a & b\\
c & d
\end{smallmatrix}\right) =a\oplus d \).

Then there exists an isomorphism
\[
\phi :(A_{1}\oplus A_{2},E_{1}\oplus E_{2})*_{D\oplus D} (M_{2}(D),F)\cong ((A_{1},E_{1})*_{D}(A_{2},E_{2}))\otimes M_{2}(\mathbb {C}),\]
so that \( \phi  \) intertwines the automorphisms \( (\alpha _{1}\oplus \alpha _{2})*(\id \otimes \alpha_1|_D)  \)
and \( (\alpha _{1}*\alpha _{2})\otimes \id \).
\end{lem}
\begin{proof}
Consider in \( M_2(D)\subseteq(A_{1}\oplus A_{2},E_{1}\oplus E_{2})*_{D\oplus D}(M_{2}(D),F) \)
the unitary \( w=\left( \begin{smallmatrix}
0 & 1\\
1 & 0
\end{smallmatrix}\right)  \). Consider the subalgebras
\begin{eqnarray*}
A_{1}' & = & a\oplus 0,\qquad a\in A_{1}\\
A_{2}' & = & w(0 \oplus a)w^{*},\qquad a\in A_{2}.
\end{eqnarray*}
Since \( w \) is free from \( A_{1}\oplus A_{2} \) with amalgamation over
\( D\oplus D \), it is easily seen from the definition of freeness and
Lemma \ref{Lemma:freenessconjugation} that the
algebras \( A_{1}' \) and \( A_{2}' \) are also free with amalgamation over
\( D\oplus 0 \). Furthermore, the restriction of \( (E_{1}\oplus E_{2})*F \)
to \( A_{i}' \) is \( E_{i} \), and hence the GNS representation associated
to each \( E_{i} \) is faithful. It follows from the embedding result
(see \S\ref{subsec:fp})
that \( A_{1}' \) and \( A_{2}' \) together generate the reduced free product
\( A_{1}*_{D}A_{2} \). It is easily seen that the algebra \( (A_{1}\oplus A_{2})*_{D\oplus D}(M_{2}(D)) \)
is generated by \( A_{1}',A_{2}' \) and \( C^{*}(1_{D}\oplus 1_{D},w)\cong M_{2}(\mathbb {C}) \),
and is isomorphic to \( (A_{1}*_{D}A_{2})\otimes M_{2}(\mathbb {C}) \) via
the map 
\begin{eqnarray*}
\phi: A_1'\ni a_1\oplus 0 & \mapsto & \left(\begin{matrix} a_1 & 0 \\ 0 & 0 \end{matrix} \right)\in A_1 \otimes M_2(\mathbb{C}),\quad a_1 \in A_1\\
\phi: A_2'\ni w 0 \oplus a_2 w^*   & \mapsto &  \left(\begin{matrix} a_2 & 0 \\ 0 & 0 \end{matrix}\right)\in A_2 \otimes M_2(\mathbb{C}),\quad a_2 \in A_2\\
\phi (w) & = & \left( \begin{array}{cc}
0 & 1\\
1 & 0
\end{array}\right) .
\end{eqnarray*}
It is clear that \( \phi  \) intertwines the automorphisms \( (\alpha _{1}\oplus \alpha _{2})*\id  \)
and \( (\alpha _{1}*\alpha _{2})\otimes 1 \).
\end{proof}

In the following lemma, $\mathcal{O}_2$ will be the Cuntz algebra~\cite{cuntz}
generated by isometries
$S_1$ and $S_2$ such that $S_1S_1^*+S_2S_2^*=1$ and $\sigma$ will denote the state on
$\mathcal{O}_2$ satisfying
\[
\sigma(S_{i_1}\cdots S_{i_p}S_{j_q}^*\cdots S_{j_1}^*)=
\begin{cases}
2^{-p}&\mbox{if }p=q\mbox{ and }i_1=j_1,\,\cdots,i_p=j_p \\
0&\mbox{ otherwise.}
\end{cases}
\]
Note that $\sigma$ is faithful.

\begin{lem}
\label{lemma:ExistenceofPartialIsometry}Let \( \sigma  \) be the state on
\( \mathcal{O}_{2} \) as above.
Let $C$ be a commutative $C^*$--algebra having a faithful state $\rho$ and a unitary $u\in C$
such that $\rho(u)=0$.
Consider the algebra \(
B=((\mathbb {C}\oplus \mathbb {C})\otimes \mathcal{O}_{2},\theta
\otimes \sigma )*(C,\rho ) \), where \(
\theta (a\oplus b)=\frac{1}{2}(a+b) \).
Then $B$ is simple and purely infinite.
Moreover, denoting by \(
E=E^{B}_{\mathbb (\mathbb{C}\oplus \mathbb {C})\otimes \mathcal{O}_{2}}
\) the conditional expectation from \( B \) onto \( \mathbb
(\mathbb{C}\oplus \mathbb {C})\otimes \mathcal{O}_{2} \) arising from the
free product construction, there exists a subalgebra \( M\subset
B \), so that:
\renewcommand{\labelenumi}{(\roman{enumi})}
\begin{enumerate}
\item  \( M\cong M_{2}(\mathbb {C}) \) in such a
  way that the element \( (a\oplus b)\otimes 1 \) corresponds to the
  diagonal matrix \( \left( \smallmatrix a & 0\\ 0 & b
  \endsmallmatrix\right)  \),
\item  \( E(M)\subset (\mathbb {C}\oplus \mathbb {C})\otimes 1 \).
\end{enumerate}
\end{lem}
\begin{proof}
Since \( \mathcal{O}_{2} \)
has trivial \( K \)--theory, it follows from Germain's exact sequence for free
products (see \S\ref{sec:Ktheoryfreeprod}) that \( K_{0}(B) \) is zero.
Once $B$ is known to be simple and purely infinite, it will follow from
Cuntz's fundamental results
\cite{cuntz:Ktheorypurelyinfinite}, that there exists a partial isometry \( w\in B \) so
that \( ww^{*}=(1\oplus 0)\otimes 1 \) and \( w^{*}w=(0\oplus 1)\otimes 1 \).
Set \( M=C^{*}(w) \). Then \( M\cong M_{2}(\mathbb {C}) \), in such a way
that \( w \) corresponds to the partial isometry \( \left(\smallmatrix
0 & 1\\
0 & 0
\endsmallmatrix\right)  \), and therefore (i) is satisfied. To see that (ii) is satisfied, note that \( M \)
is the linear span of \( w,w^{*},ww^{*},w^{*}w \); hence it is sufficient to
check that \( E(w),E(w^{*}w),E(ww^{*}) \) and \( E(w^{*})=E(w)^{*} \) all
lie in \( (\mathbb {C}\oplus \mathbb {C})\otimes 1 \). This is clearly true
of \( E(w^{*}w) \) and \( E(ww^{*}) \), since \( ww^{*} \) and \( w^{*}w \)
lie in \( (\mathbb {C}\oplus \mathbb {C})\otimes 1 \). Furthermore, since \( (ww^{*})E(w)=E(ww^{*}w)=E(w)w^{*}w) \),
and the projections \( ww^{*} \) and \( w^{*}w \) are orthogonal, it follows
that \( ww^{*}E(w)ww^{*} \) and \( w^{*}wE(w)w^{*}w \) are both zero. Since
\( ww^{*} \) and \( w^{*}w \) lie in the center of \( (\mathbb {C}\oplus \mathbb {C})\otimes \mathcal{O}_{2} \)
and \( ww^{*}+w^{*}w=1 \), it follows that \( E(w)=0 \), and hence \( E(w^{*})=0 \),
so that (ii) is satisfied as well.

Note that \( B\cong (\mathcal{O}_{2}\oplus \mathcal{O}_{2},\sigma ')*(C,\rho ) \),
where \( \sigma '=\frac{1}{2}(\sigma +\sigma ) \).
We shall now apply Theorem 3.1 of~\cite{dykema:piII} to show that $B$ is simple and purely infinite;
the following five observations show that the hypotheses of this theorem are satisfied.

\renewcommand{\labelenumi}{\arabic{enumi}.}
\begin{enumerate}

\item $\sigma'$ and $\rho$ are faithful states.

\item Let $v=S_2S_2S_1^*\oplus0 \in \mathcal{O}_2 \oplus\mathcal{O}_2$;
  then $v$ is a partial isometry belonging to the spectral subspace of
  the state $\sigma'* \rho$ corresponding to $2$, i.e.\ $F(xv)=2F(vx)$ for all $x\in B$.

\item Let $q_1=v^*v$ and $q_2=vv^*$; then $q_1\perp q_2$ and we have $F(q_1)=\frac{1}{4}$, $F(q_2)=\frac{1}{8}$.

\item Consider
\begin{align*}
B_1&=C^*(q_1,q_2,u)\subseteq B \\
B_2&=C^*(uq_1u^*,q_2)\subseteq B_1.
\end{align*}

By Lemma~\ref{Lemma:freenessconjugation}, $uq_1u^*$ and $q_2$ are free with respect to $F$, so
\[
B_2\cong(\underset{3/4}{\mathbb{C}}\oplus\overset{uq_1u^*}{\underset{1/4}{\mathbb{C}}})
*(\underset{7/8}{\mathbb{C}}\oplus\overset{q_2}{\underset{1/4}{\mathbb{C}}})
\]
and it is well known (cf. Proposition 2.7 of~\cite{dykema:simplicityfreeprod}) that
\[
B_2\cong\mathbb{C}\oplus C([0,1],M_2(\mathbb{C}))\oplus\mathbb{C}
\]
with
\begin{align*}
uq_1u^*&\sim1\oplus\left(\begin{array}{ccc}1&0\\0&0\end{array}\right)\oplus0 \\
q_2&\sim0\oplus\left(\smallmatrix t&\sqrt{t(1-t)}\\ \sqrt{t(1-t)}&1-t\endsmallmatrix\right)\oplus0.
\end{align*}
Therefore, $q_2$ is equivalent in $B_1$ to a subprojection of $q_1$ and $q_2B_1q_2$ contains the diffuse abelian subalgebra $q_2B_2q_2$.

\item The centralizer of $\sigma'$ contains a diffuse abelian subalgebra; hence by Proposition 3.2 of~\cite{dykema:simplicityfreeprod}, $B$ is simple and thus $q_1+q_2$ is full in $B$.
\end{enumerate}

The above facts allow us to apply Theorem 3.1 of~\cite{dykema:piII}, and we conclude that $B$ is simple and purely infinite.
\end{proof}

\begin{lem}
\label{lemma:XYZZY}Let $W$ be a unital $C^*$--algebra and
let \( E^{X}_{W}:X\to W \) and \( E^{Y}_{W}:Y\to W \)
be \( W \)--probability spaces, with automorphisms \( \alpha _{X} \) and \( \alpha _{Y} \),
respectively, which agree on \( W \). Assume that \( W\subset Z\subset X \)
is an \( \alpha _{X} \)--invariant subalgebra, and let \( E^{Z}_{W} \) be the
restriction of \( E^{X}_{W} \) to \( Z \). Assume that there exists an \( E^{X}_{W} \)--preserving
\( \alpha _{X} \)--invariant conditional expectation \( E^{X}_{Z}:X\to Z \),
and assume that the GNS representations associated to \( E^{X}_{W},E^{Z}_{W},E^{X}_{Z} \)
and \( E^{Y}_{W} \) are faithful. Then there exists an isomorphism of the reduced
free products
\[
X*_{Z}(Z*_{W}Y)\cong X*_{W}Y,\]
which intertwines the automorphisms \( \alpha _{X}*(\alpha _{X}|_{Z}*\alpha _{Y}) \)
and \( \alpha _{X}*\alpha _{Y} \).
\end{lem}

\begin{proof}
Consider the free product conditional expectation \(E_Z : (X*_Z (Z*_W
Y) \to Z\).  The GNS representation associated to \( E_Z \) is
faithful, by definition.  Let \( E_W = E^Z_W \circ E_Z \).  By Lemma
\ref{lem:condexp}, \( E_W \) also gives rise to a faithful GNS representation.  
By Lemma 2.6 in \cite{shlyakht:cpentropy}, \( X \) and \( Y \) are
free with respect to \( E_W \).  It follows from the assumptions and
the embedding result (see \S\ref{subsec:fp})
that the \(C^*\)--algebra generated by \(X\) and \(Y\) in the GNS
representation of \(X*_{Z}(Z*_{W}Y) \) associated to \( E_W \) is
isomorphic to \( X*_W Z\).  The desired isomorphism follows.

\end{proof}

\begin{lem} \label{lem:mainembedding}
Let \( E_{i}:A_{i}\to D \), \( i=1,2 \) be \( D \)--probability spaces, with
automorphisms \( \alpha _{i} \). Assume that the GNS representations associated
to \( E_{i} \) are faithful, and that \( \alpha _{1}|_{D}=\alpha _{2}|_{D} \).

Let \( E:A_{1}\oplus A_{2}\to D=(d\oplus d:d\in D)\subset A_{1}\oplus A_{2} \)
be given by \[ E(a_{1}\oplus a_{2})=\frac{1}{2}([E(a_{1})+E(a_{2})]\oplus [E(a_{1})\oplus E(a_{2})]). \]
Let \( \sigma :\mathcal{O}_{2}\to \mathbb {C} \)
be as in Lemma~\ref{lemma:ExistenceofPartialIsometry}.
Let \( C \) be any \( C^{*} \)--algebra
with a state \( \hat{\rho } \), giving rise to a faithful GNS representation,
and containing a unitary $u$, s.t. \(\rho=\hat\rho|_{C^*(u)}\) is
faithful and \( \rho(u)=0 \).

Then there exists an embedding
\begin{eqnarray*}
\lambda :(A_{1},E_{1})*_{D}(A_{2},E_{2})\to B=((A_{1}\oplus A_{2})\otimes \mathcal{O}_{2}),E\otimes \sigma )*_{D}(C\otimes _{\min }D,\hat{\rho }\otimes \id ), & 
\end{eqnarray*}
so that \( \lambda  \) intertwines \( (\alpha _{1}*\alpha _{2}) \) and \( ((\alpha _{1}\oplus \alpha _{2})\otimes \id )*(\id \otimes \alpha _{1}|_{D}) \).
\end{lem}

\begin{proof}
Using the embedding result (see \S\ref{subsec:fp}), we can reduce
to the case that \( C=C^*(u) \) and \( \hat{\rho }=\rho  \), by replacing
\( B \) with the algebra generated by \( (A_{1}\oplus A_{2})\otimes \mathcal{O}_{2} \)
and \( C^*(u) \otimes D\subset C\otimes _{\min }D \). 

We have the following sequence of covariant inclusions, justified below:
\begin{eqnarray*}
A_{1}*_{D}A_{2} & \hookrightarrow  & M_{2}\otimes (A_{1}*_{D}A_{2})\\
 & \stackrel{\textrm{Lemma \ref{lemma:embedfreeprodintodirectsum}}}{\cong } & M_{2}\otimes D*_{D\oplus D}(A_{1}\oplus A_{2})\\
 & \stackrel{(\textrm{a})}{\hookrightarrow } & \left( [(\mathbb {C}^{2})\otimes \mathcal{O}_{2})*C]\otimes D\right) *_{(\mathbb {C}^{2})\otimes \mathcal{O}_{2}\otimes D}[(A_{1}\oplus A_{2})\otimes \mathcal{O}_{2}]\\
 & \cong  & \left( (\mathbb {C}^{2}\otimes \mathcal{O}_{2}\otimes D)*_{D}(C\otimes D)\right) *_{(\mathbb {C}^{2})\otimes \mathcal{O}_{2}\otimes D}[(A_{1}\oplus A_{2})\otimes \mathcal{O}_{2}]\\
 & \stackrel{(\textrm{b})}{\cong } & C\otimes D*_{D}[(A_{1}\oplus A_{2})\otimes \mathcal{O}_{2}].
\end{eqnarray*}
Inclusion (a) is implied by Lemma \ref{lemma:ExistenceofPartialIsometry}, together with
the embedding result (see \S\ref{subsec:fp}). 

Isomorphism (b) is implied by Lemma \ref{lemma:XYZZY}.

\end{proof}

\begin{thm}
\label{Theorem:mainresult}
Let $D$ be a finite dimensional C$^*$--algebra and let \( \alpha _{j} \), be an automorphism
of a \( D \)--probability space \( (A_{j},E_{j}:A_{j}\to D) \), ($j\in J$,
$J$ a finite or countably infinite set).
Assume that
the GNS representations associated to \( E_{j} \) ($j\in J$) are faithful.
Assume that \( \alpha _{j}|_{D}=\alpha _{i}|_{D} \) for all \( i,j\in J \).
Let $(A,E)=*_{D} \left((A_{j},E_{j}),j\in J \right)$ and let \(*_j \alpha _{j} \)
denote the free product automorphism of $A$.
Then
\begin{equation}
\label{eq:HTfp}
\HT (*_j \alpha _{j})=\sup_{j\in J} (\HT (\alpha _{j})).
\end{equation}
\end{thm}
\begin{proof}
Because of the embedding result (see~\S\ref{subsec:fp}) and the behavior of entropy with 
respect to inductive limits, it suffices to prove the statement when
$J=\{1,2\}$
is a set with two elements.

By monotonicity of $\HT$, the inequality $\ge$ in~\eqref{eq:HTfp} is clear.
Let $\mathcal{T}$ be the Toeplitz algebra generated by the unilateral shift $\ell$
and denote by $\hat\rho$ the vacuum expectation on $\mathcal T$.
Voiculescu~\cite{dvv:symmetries} showed that $\ell^*+\ell$ is a semicircular element with its spectrum an interval.  Hence
by functional calculus $C=C^*(\ell^*+\ell)$ contains a unitary $u$
such that $\hat\rho(u^k)=0$, for all $k\neq 0$.
It is not difficult to see that $C$ contains no nonzero compact
operator, and thus $\rho=\hat{\rho}|_{C}$
is faithful.
By Lemma~\ref{lem:mainembedding} and Lemma~\ref{lemma:embedfreeprodintodirectsum}, the amalgamated free product dynamical
system \( (A_{1}*_{D}A_{2},\alpha _{1}*\alpha _{2}) \) can be covariantly embedded into $(B,\beta)$, where
\[
(B,F)=((A_1\oplus A_2)\otimes\mathcal{O}_2,\widetilde{E}\otimes\sigma)*_D(\mathcal{T}\otimes D,\hat{\rho}\otimes\id_D)
\]
and where
\[
\beta=((\alpha_1\oplus\alpha_2)\otimes\id_{\mathcal{O}_2})*(\id_{\mathcal T}\otimes(\alpha_1|_D)).
\]
So by monotonicity again, $\HT(\alpha_1*\alpha_2)\le\HT(\beta)$.
Using Theorem~\ref{Theorem:freeprodToeplitz} we have
\[
\HT(\beta)=\HT((\alpha_1\oplus\alpha_2)\otimes\id_{\mathcal{O}_2})=\HT(\alpha_1\oplus\alpha_2)
=\max(\HT(\alpha_{1}),\HT(\alpha_{2})).
\]
\end{proof}

\begin{rems}{}\hspace{1ex}
\label{rmq:afalgebra}

\medskip

\noindent
\thethm.1.
The only reason that the hypotheses of Theorem \ref{Theorem:mainresult} require
finite dimensionality of \( D \) is because the proof appeals to Theorem \ref{Theorem:freeprodToeplitz}.
If one could prove a version of Theorem \ref{Theorem:freeprodToeplitz} for
more general types of \( D \), then Theorem \ref{Theorem:mainresult} would
be valid for those more general algebras \( D \), with the same proof.
\medskip

\noindent
\thethm.2.
Theorem~\ref{Theorem:mainresult} also holds for amalgamated free products of
injective {\em endomorphisms} of non--commutative probability spaces. 
This can be seen by realizing an injective endomorphism as the restriction
of an automorphism (see e.g. \cite{shlyakht-dykema:exactness}), 
and then utilizing our result for automorphisms.
\medskip

\noindent
\thethm.3.
One can actually prove a partial version of Theorem~\ref{Theorem:mainresult} with amalgamation
taking place over an AF algebra.  Assume that $E_1:A_1\to D$ and $E_2:A_2\to D$ are $D$--probability spaces, with automorphisms $\alpha_{1}$, $\alpha_2$,
which restrict to the same automorphism on $D$.  Assume that $D=\overline{\cup_i D^{(i)}}$, with $D^{(i)}\subset D^{(i+1)}$ finite dimensional. 
   Assume 
further that $A_j =\overline{\cup_i A_j^{(i)}}$, $A_j^{(i)}\subset A_j^{(i+1)}$, and that $E_j (A_j^{(i)}) = D^{(i)}$ and $D^{(i)}\subset A_j^{(i)}$, $j=1,2$.
Let $\alpha_1 * \alpha_2$ denote the free product automorphism on $A_1
*_D A_2$. Then $\HT(\alpha_1 *
\alpha_2)=\max(\HT(\alpha_1),\HT(\alpha_2))$. 
This is due to the fact that under these assumptions, \( A_1 *_D A_2 \) is the direct limit of
\( A^{(i)}_1 *_{D^{(i)}} A^{(i)}_2 \).
\end{rems}

\section{Applications}

\subsection{Shifts on infinite free products and twisted free permutations.}

Let \( A \) be a unital \( C^{*} \)--algebra, \( D \) a finite dimensional
unital subalgebra and let \( E:A\to D \) be a conditional expectation with a faithful
GNS representation.
Let $I$ be a nonempty set and for every $i\in I$ let \( (A_{i},E_i) \), 
be a copy of \( (A,E) \).
Let
\[
(B,F)=(*_{D})_{i\in I}(A_{i},E_i)
\]
be the reduced amalgamated free product.
Let \( \sigma :I\to I \)
be a bijection.
The {\em free permutation}
arising from $\sigma$ is the automorphism $\beta$ of $B$ that permutes the copies
of $A$ inside $B$ by sending
\( A_i \) to \( A_{\sigma(i)} \).
(The {\em free shift} is the free permutation arising from the shift on $I=\mathbb Z$.)
The following theorem is a generalization of results
of St\o{}rmer \cite{stormer:freeshiftstypeII}, Brown and Choda \cite{brown-choda:entropycrossed}
and Dykema \cite{dykema:entropyfreeprod} and a partial generalization of other results
of St\o{}rmer, \cite{stormer:freeshifts}.
It answers affirmatively in the case
of amalgamation over a finite dimensional \( C^{*} \) algebra Question~11
of~\cite{dykema:entropyfreeprod}:

\begin{subthm}
\label{theorem:freeshiftzero}
If $A$ is an exact unital C$^*$--algebra and $D$ is finite dimensional
then for every free permutation $\beta$ of $B=(*_D)_{i\in I}A$ one has
\( \HT (\beta )=0 \).
\end{subthm}
\begin{proof}
Consider the reduced group $C^*$--algebra $C^*_r(\mathbb F_{|I|})$ of the nonabelian free
group on $|I|$ generators, where $|I|$ is the cardinality of $I$.
Let $(u_i)_{i\in I}$ be the unitary generators of $C^*_r(\mathbb F_{|I|})$ corresponding
to the free generators of $\mathbb F_{|I|}$.
Let $\sigma_*\in\Aut(C^*_r(\mathbb F_{|I|}))$ be the automorphism such that $\sigma_*(u_i)=u_{\sigma(i)}$.
Then from~\cite{brown-choda:entropycrossed} and~\cite{dykema:entropyfreeprod}, $\HT(\sigma_*)=0$,
Let
\[
(\widetilde B,\widetilde F)=(C^*_r(\mathbb F_{|I|})\otimes D,\tau\otimes\id_D)*_D(A,E),
\]
where $\tau$ is the canonical tracial state on $C^*_r(\mathbb F_{|I|})$, and let
\[
\tilde\beta=(\sigma_*\otimes\id_D)*\id_A\in\Aut(\widetilde B).
\]
Then by Theorem~\ref{Theorem:mainresult}, $\HT(\tilde\beta)=\HT(\sigma_*)=0$.
Let $B'=C^*(\bigcup_{i\in I}u_iAu_i^*)\subseteq\widetilde B$.
Then $\widetilde F(u_iau_i^*)=E(a)$ for all $a\in A$ and $i\in I$ and,
by Lemma~\ref{Lemma:freenessconjugation}, the family $(u_iAu_i^*)_{i\in I}$
is free with respect to $\widetilde F$.
It is not difficult to see, though somewhat tedious to write down in detail, 
that the inclusion representation of $B'$ on $L^2(\widetilde B,\widetilde F)$
is a multiple of the GNS representation, $\rho$, of $B'$ on $L^2(B',\widetilde F|_{B'})$.
Indeed, one chooses vectors $\xi_i\in L^2(\widetilde B,\widetilde F)$ such that
$\bigcup_iB'\xi_i$ spans a dense subset of $L^2(\widetilde B,\widetilde F)$ and so that
the representation of $B'$ on $B'\xi_i$ is equivalent to $\rho$, for all $i$.
Hence there is an isomorphism $\pi:B\to B'$ sending the copy of $A_i$ in $B$ to $u_iAu_i^*$.
We see that the automorphisms $\beta$ and $\tilde\beta|_{B'}$ are conjugate via $\pi$.
Hence, by monotonicity of $\HT$ we have
$\HT(\beta)\le\HT(\tilde\beta)=0$.
\end{proof}

\begin{subrem}
Just as in \ref{rmq:afalgebra}.3, it is possible to weaken the
hypothesis of the preceding Theorem by requiring that $(A,D,E:A\to D)$
is an inductive limit of $(A^{(i)},D^{(i)},E_i)$, with $D^{(i)}$
finite dimensional.
\end{subrem} 

\begin{subdefn}
Let $D$ be a unital C$^*$--algebra, let $I$ be a set and for every $i\in I$
let $(A_i,E_i)$ be a $D$--probability space such that the GNS representation
associated to $E_i$ is faithful.
Let
\[
(A,E)=(*_D)_{i\in I}(A_i,E_i)
\]
be the reduced amalgamated free product.
Let $\sigma:I\to I$ be a bijection such that for every $i\in I$
there is an isomorphism $\alpha_i:A_i\to A_{\sigma(i)}$
such that $\alpha_i(D)=D$ and $E_{\sigma(i)}\circ\alpha_i=E_i$.
Assume that $\alpha_i|_D$ is the same for all $i$.
Then there is an automorphism $\alpha$ of $A$ sending the copy of $A_i$
in $A$ onto the copy of $A_{\sigma(i)}$ in $A$ via $\alpha_i$.
An automorphism $\alpha$ arising in this way is called a {\em twisted free
permutation} of $A$.
\end{subdefn}

The next theorem generalizes Theorem~\ref{theorem:freeshiftzero}
and also some results of~\cite{dykema:entropyfreeprod}.
\begin{subthm}
Suppose that $D$ is finite dimensional and $\alpha\in\Aut(A)$ is a twisted free
permutation as described above.
If the permutation $\sigma$ of $I$ has no cycles then $\HT(\alpha)=0$.
Otherwise, whenever $c$ is a cycle of $\sigma$ let $\ell(c)$ be the length of the cycle,
let $i\in I$ be one of the elements moved by the cycle and let $\beta_c\in\Aut(A_i)$
be the restriction of $\alpha^{\ell(c)}$ to the copy of $A_i$ in $A$.
(Note that $\beta_c$ depends on $i$ only up to conjugation.)
Then
\begin{equation}
\label{eq:HTtfperm}
\HT(\alpha)=\sup_c\frac{\HT(\beta_c)}{\ell(c)},
\end{equation}
where the supremum is over all cycles $c$ of $\sigma$.
\end{subthm}
\begin{proof}
If $\sigma$ has no cycles then $\alpha$ is (conjugate to) a free permutation so $\HT(\alpha)=0$
by Theorem~\ref{theorem:freeshiftzero}.
In general, by making a cycle decomposition of $\sigma$ and using Theorem~\ref{Theorem:mainresult},
we see that in order
to prove~\eqref{eq:HTtfperm} we may without loss of generality assume that $\sigma$ itself is a
cyclic permutation, $\sigma=c$, of a finite set $I$.
However, then $\alpha^{\ell(c)}=*_{i\in I}\gamma_i$, where each $\gamma_i\in\Aut(A_i)$ is conjugate to
$\beta_c$.
Hence again applying Theorem~\ref{Theorem:mainresult}
we have $\HT(\alpha)=\HT(\alpha^{|I|})/|I|=\HT(\beta_c)/\ell(c)$.
\end{proof}

\subsection{The CNT variational principle.}

In classical ergodic theory an important result connecting topological and measurable
entropy is the variational principle. 

\begin{subdefn}
Let \( (A,\alpha ) \) be a unital exact \( C^{*} \)--dynamical system.  We
say that \( (A,\alpha ) \) satisfies a CNT--variational principle if
\[
\HT (\alpha )=\sup _{\phi }h_{\phi }(\alpha ),\]
where the supremum is taken over all \( \alpha  \)--invariant states on \( A \)
and \( h_{\phi }(\alpha ) \) denotes the CNT--entropy of \( \alpha  \) with
respect to \( \phi  \) (cf. \cite{connes-narnhofer-thirring:entropy}).
\end{subdefn}
By \cite[Prop. 9]{dykema:entropyfreeprod} we always have the inequality
\[
\HT (\alpha )\geq \sup _{\phi }h_{\phi }(\alpha ).\]

Not every \( C^{*} \)--dynamical system satisfies a CNT--variational principle.
 In \cite{narnhofer-stormer-thirring:entropytensorfails} an example of a highly
non--asymptotically abelian system was given for which \( \HT (\alpha )\geq (\log 2)/2 \)
 while \( h_{\phi }(\alpha )=0 \) for the unique \( \alpha  \)--invariant (tracial)
state.  But, since \( \HT (\cdot ) \) agrees with classical topological entropy
and \( h_{\phi }(\cdot ) \) agrees with classical Kolmogorov--Sinai entropy
when \( A \) is abelian, the classical variational principle says that if \( A \)
is unital and abelian then \( (A,\alpha ) \) always satisfies a CNT--variational
principle.  The list of examples of noncommutative dynamical systems which satisfy
a CNT--variational principle is also rapidly growing (cf. \cite[4.7]{dvv:topentropy},
\cite[4.6, 4.7]{choda:entropyCuntzEndomorphism}, \cite[3.6, 3.7]{brown:entropy},
\cite{boca-goldstein:entropyCuntzPimsner}, \cite{pinzari-watatani-yonetani:variationalprincipleentropy},
\cite{stormer-neshveev:variationalprincipleentropy}).  Moreover, the class
of \( C^{*} \)--dynamical systems which satisfy a CNT--variational principle
is closed under taking (minimal) tensor products (cf. \cite[Lemma 3.4]{dvv-stormer:entropyCAR})
and crossed products by \( {\mathbb Z} \) (cf. \cite[Theorem 3.5]{brown:entropy}).
Unfortunately it is \emph{not} closed under taking quotients or subalgebras
(even those with a conditional expectation onto them --- simply take a direct
sum of something abelian with large entropy and the example from \cite{narnhofer-stormer-thirring:entropytensorfails})
and it not yet known what happens in extensions.  However we now show that it
is also closed under taking reduced free products.  The next theorem also gives
lots of examples of non--asymptotically abelian dynamical systems for which the
CNT--variational principle holds.   

\begin{subthm}
Let $E_j : A_j \to D$, $j=1,2$ be non--commutative probability spaces
with automorphisms $\alpha_1$ and $\alpha_2$.  Assume that $D$ is finite dimensional
and that $\alpha_1|_D=\alpha_2|_D$, and that $E_j$ give rise to faithful
GNS representations.  If both $(A_1,\alpha_1)$ and $(A_2,\alpha_2)$ satisfy
the CNT variational principle, then so does $(A_1 *_D A_2,\alpha_1 *\alpha_2)$.
\end{subthm}
\begin{proof}
Assume without loss of generality that \( \HT (\alpha_1 )\geq \HT (\beta_1 ) \).
 Given \( \epsilon >0 \) we can find an \( \alpha_1  \)--invariant state \( \gamma \in S(A_1) \)
such that \( \HT (\alpha_1 )\leq h_{\gamma }(\alpha_1 )+\epsilon  \).  

Let \( E:A_1*A_2\to A_1 \) be the canonical conditional expectation
and define a state \( \tilde{\gamma }=\gamma \circ E\in S(A_1 *_D A_2) \).
 Then one checks that \( \tilde{\gamma }\circ (\alpha_1 *\alpha_2 )=\tilde{\gamma } \)
and \( \tilde{\gamma }\circ E=\tilde{\gamma } \).  But under these conditions,
CNT--entropy is also monotone (cf. \cite[III.6]{connes-narnhofer-thirring:entropy})
and so \( \HT (\alpha_1 *\alpha_2 )=\HT (\alpha_1 )\leq h_{\gamma }(\alpha_1 )+\epsilon \leq h_{\tilde{\gamma }}(\alpha_1 *\alpha_2 )+\epsilon  \).
 Since \( \epsilon  \) was arbitrary, we are done.  
\end{proof}
A particularly interesting class of dynamical systems for which one would like
to have a CNT--variational principle is those arising from automorphisms of discrete groups.  For
any discrete group \( G \) let \( C^{*}_{r}(G) \) be the reduced group \( C^{*} \)--algebra
and \( \tau _{G} \) the canonical trace on \( C^{*}_{r}(G) \).  \( G \) is called exact if \( C^{*}_{r}(G) \) is exact.  If \( \gamma :G\to G \)
is a group automorphism then there is an induced automorphism \( \hat{\gamma }\in \Aut (C^{*}_{r}(G)) \)
such that \( \tau _{G}\circ \hat{\gamma }=\tau _{G} \).  If \( G \) is abelian
then a classical theorem of Berg \cite{berg:maximalentropyHaarmeasure} implies
that \( \HT (\hat{\gamma })=h_{\tau _{G}}(\hat{\gamma }) \).  

\begin{subthm}
Let $H$ be a finite group and let
$G_1$ and $G_2$ be discrete exact groups having a $H$ as a common subgroup.
Suppose that $\gamma_i\in\Aut(G_i)$, ($i=1,2$) are automorphisms preserving $H$
such that $\gamma_1|_H=\gamma_2|_H$.
Let $\gamma_1*\gamma_2$ denote the resulting automorphism of the free product $G_1*_HG_2$
with amalgamation over $H$.
If $\HT(\hat\gamma_i)=h_{\tau_{G_i}}(\hat\gamma_i)$, ($i=1,2$), then
$\HT((\gamma_1*\gamma_2)\widehat{\;})=h_{\tau_{G_1*_HG_2}}((\gamma_1*\gamma_2)\widehat{\;})$.
\end{subthm}
\begin{subcor}
If \( G_{1},G_{2} \) are discrete exact groups with automorphisms \( \gamma _{i}\in \Aut (G_{i}) \),
(\( i=1,2 \)), and if \( \HT (\hat{\gamma }_{i})=h_{\tau _{G_{i}}}(\hat{\gamma }_{i})\) ,(\(i=1,2 \)),
then \( \HT (\hat{\gamma }_{1}*\hat{\gamma }_{2})=h_{\tau _{G_{1}*G_{2}}}(\hat{\gamma }_{1}*\hat{\gamma }_{2}) \).
\end{subcor}
It would be interesting to know whether or not the above theorem or its corollary can be extended
to the dual entropy defined in \cite{brown-germain:dualentropy} as well.

\subsection{Entropy--preserving embeddings.}
\label{sec:htpresemb}

Kirchberg first proved that every separable exact \( C^{*} \)--algebra is isomorphic
to a subalgebra of the Cuntz algebra on two generators (cf. \cite{kirchberg-phillips:embeddingCuntz}).
 In \cite[Remark 2.3]{brown:entropy} it was asked whether or not one can always
find a unital embedding \( \rho :A\rtimes _{\alpha }{\mathbb Z}\hookrightarrow \mathcal{O}_{2} \)
such that \( \HT _{\mathcal{O}_{2}}({\Ad }\rho (u))=\HT (\alpha ) \), where
\( u\in A\rtimes _{\alpha }{\mathbb Z} \) is the implementing unitary.  We
now solve this problem affirmatively in the case that \( A \) is nuclear and
there exists an \( \alpha  \)--invariant state \( \phi \in S(A) \) with faithful
GNS representation. We also show that if \( (A,\alpha ) \) is \emph{any} nuclear
\( C^{*} \)--dynamical system then there always exists an entropy preserving
covariant embedding into the Cuntz algebra on infinitely many generators.  Since
many \( C^{*} \)--algebras are stable under tensor products with \( \mathcal{O}_{\infty } \)
we thus get entropy preserving embeddings into a very large class of \( C^{*} \)--algebras.
 It follows that the topological entropy invariant of all such algebras is \( [0,\infty ] \).

We begin with a simple proposition. 

\begin{subprop}
\label{prop:innerfreeeautomorphism}Assume \( A \) is a unital exact \( C^* \)--algebra,
\( \alpha \in \Aut(A) \) and there exists an \( \alpha  \)--invariant state
\( \phi \in S(A) \) with faithful GNS representation.  Let \( B \) be unital
and exact, let \( \psi \in S(B) \) have faithful GNS representation, let \( E:A\rtimes _{\alpha }{\mathbb Z}\to A \)
be the canonical conditional expectation and let \( u\in A\rtimes _{\alpha }{\mathbb Z} \)
be the implementing unitary. Regarding \( u \) as a unitary in \( (A\rtimes _{\alpha }{\mathbb Z},\phi \circ E)*(B,\psi ) \)
we have \( \HT ({\Ad }u)=\HT (\alpha ) \). 
\end{subprop}
\begin{proof}
Since \( E:A\rtimes _{\alpha }{\mathbb Z}\to A \) is faithful, by Lemma~\ref{lem:condexp} the GNS representation
of \( \phi \circ E \) is faithful. 

Consider the \( C^{*} \)--algebra \( (A,\phi )*\big (*_{\mathbb Z}(B,\psi )\big )\rtimes _{\alpha *S}{\mathbb Z} \),
where \( S:*_{\mathbb Z}(B,\psi )\to *_{\mathbb Z}(B,\psi ) \) is the free
shift and let \( v\in (A,\phi )*\big (*_{\mathbb Z}(B,\psi )\big )\rtimes _{\alpha *S}{\mathbb Z} \)
be the implementing unitary.  It is fairly easy to see (cf. the proof of \cite[Claims 2 and 4]{choda-dykema:pisumfromfreeIII})
that there exists a \( * \)--isomorphism
\[
\rho :(A\rtimes _{\alpha }{\mathbb Z},\phi \circ E)*(B,\psi )\to
(A,\phi )*\big(\big (*_{\mathbb Z}(B,\psi )\big )\rtimes _{\alpha *S}{\mathbb Z}\big),
\]
such that \( \rho (u)=v \).  Hence \( \HT ({\Ad }u)=\HT ({\Ad }v) \).  But
\cite[Theorem 3.5]{brown:entropy}, Theorem \ref{theorem:freeshiftzero} above
and our main result imply that \( \HT ({\Ad }v)=\HT (\alpha *S)=\HT (\alpha ) \).  \end{proof}

\begin{subthm}
Let \( A \) be a unital separable nuclear \( C^{*} \)--algebra and \( \alpha \in \Aut (A) \)
be an automorphism such that there exists an \( \alpha  \)--invariant state
\( \phi \in S(A) \) with faithful GNS representation.  If \( \mathcal{O}_{2} \)
denotes the Cuntz algebra on two generators then there exists a unital embedding
\( \rho :A\rtimes _{\alpha }{\mathbb Z}\hookrightarrow {\mathcal{O}}_{2} \)
such that \( \HT _{{\mathcal{O}}_{2}}(\Ad \rho (u))=\HT (\alpha ) \), where
\( u\in A\rtimes _{\alpha }{\mathbb Z} \) is the implementing unitary. 
\end{subthm}
\begin{proof}
Replacing \( (A,\alpha ) \) with \( (A\otimes C({\mathbb T}),\alpha \otimes \id _{C({\mathbb T})}) \),
if necessary, we may assume that there exists a Haar unitary in the centralizer
of \( \phi  \). In this setting, \cite[Prop. 3.2]{dykema:simplicityfreeprod}
ensures that \( (A\rtimes _{\alpha }{\mathbb Z},\phi \circ E)*(\mathcal{T},v) \)
is a simple \( C^{*} \)--algebra, where \( \mathcal{T} \) is the Toeplitz algebra
and \( v \) is the vacuum state. Moreover, this free product is nuclear since
it is isomorphic to a Cuntz--Pimsner algebra over \( A\rtimes _{\alpha }{\mathbb Z} \).
 Thus Kirchberg's absorption theorem for \( \mathcal{O}_{2} \) (i.e. \( B\otimes \mathcal{O}_{2}\cong \mathcal{O}_{2} \),
for any simple, separable, unital, nuclear \( C^{*} \)--algebra \( B \), cf.
\cite{kirchberg-phillips:embeddingCuntz}) together with the previous proposition
implies the result. 
\end{proof}
We now turn to covariant entropy preserving embeddings into \( \mathcal{O}_{\infty } \).

\begin{subthm}
\label{theorem:cuntzalgebrahasallentropy}Let \( A \) be any separable nuclear
\( C^{*} \)--algebra and \( \alpha \in \Aut (A) \) be an automorphism.  If
\( \mathcal{O}_{\infty } \) denotes the Cuntz algebra on infinitely many generators
then there exists an automorphism \( \beta \in \Aut (\mathcal{O}_{\infty }) \)
and a nonunital \( * \)--monomorphism \( \pi :A\to \mathcal{O}_{\infty } \) such that
\( \beta \circ \pi =\pi \circ \alpha  \) and \( \HT (\beta )=\HT (\alpha ) \).
\end{subthm}
\begin{proof}
By adding a unit to \( A \) and replacing \( A \) by its crossed product by
\( \alpha  \), we may assume that \( A \) is unital and \( \alpha  \) is
inner.  We may also assume that there exists an \( \alpha  \)--invariant state
\( \phi \in S(A) \) with faithful GNS representation.  Indeed, it follows from
Theorem \ref{theorem:entropyextensionestimate} that the covariant embedding
\( (A,\alpha )\hookrightarrow (A\oplus A+\mathcal{K}(H\oplus H),{\Ad }(u\oplus 1_{A})) \),
where we regard \( A\subset B(H) \) and \( u\in A \) is a unitary which implements
\( \alpha  \), is entropy preserving.  Now any vector state arising from the
second copy of \( H \) will be \( {\Ad }(u\oplus 1_{A}) \)--invariant and have
faithful GNS representation since it is a cyclic vector. 

So, we assume that \( A \) is unital, \( \alpha  \) is inner and \( \phi \in S(A) \)
is \( \alpha  \)--invariant with faithful GNS representation. Consider the free
product
\[
B=\bigg ((A\otimes C[0,1]\otimes \mathcal{O}_{2}\oplus C[0,1],\eta )\bigg )*(\mathcal{T},v).\]
Here \( \eta  \) is the average of Lebesgue measure on the second copy of \( C[0,1] \)
and the tensor product of \( \phi  \), Lebesgue measure and an arbitrary faithful
state on \( \mathcal{O}_{2} \), and \( v \) is the vacuum state on the Toeplitz
algebra \( \mathcal{T} \).  Since there is a Haar unitary in the centralizer
of \( \eta  \) (coming from the copies of \( C[0,1] \)), \cite[Prop. 3.2]{dykema:simplicityfreeprod}
implies that \( B \) is simple.  Moreover, \( B \) is nuclear being isomorphic
to a Cuntz--Pimsner algebra over \( A\otimes C[0,1]\otimes \mathcal{O}_{2}\oplus C[0,1] \)
(cf. \cite{shlyakht-dykema:exactness}, \cite{germain:KCuntz-Pimsner}).  Since
\( \mathcal{O}_{2} \) is \( KK \)--equivalent to zero, \( A\otimes C[0,1]\otimes \mathcal{O}_{2}\oplus C[0,1] \)
satisfies the Universal Coefficient Theorem.  Since \( \mathcal{T} \) is type
I it also satisfies the UCT and hence, by a result of Germain~\cite{germain:Kcommutatorideal}
so does \( B \).  Moreover, it follows
from Germain's six term exact sequence for \( K \)--theory (see \S\ref{sec:Ktheoryfreeprod})
that \( B \) has the \( K \)--theory of \( \mathcal{O}_{\infty } \).  Hence
\( B\otimes \mathcal{O}_{\infty } \) is a simple, unital, purely infinite,
nuclear \( C^{*} \)--algebra which satisfies the UCT and has the \( K \)--theory
of \( \mathcal{O}_{\infty } \).  So, by the classification results of Kirchberg--Phillips
(\cite{kirchberg:classificationPISN}, \cite{phillips:classificationPISN}),
\( B\otimes \mathcal{O}_{\infty }\cong \mathcal{O}_{\infty } \).  By our main
result we have that the automorphism \( \gamma =(\alpha \otimes \id _{C[0,1]}\otimes \id _{\mathcal{O}_{2}}\oplus \id _{C[0,1]})*\id _{\mathcal{T}}\in \Aut (B) \)
has the same entropy as \( \alpha  \).  Hence defining \( \beta =\gamma \otimes \id _{\mathcal{O}_{\infty }}\in \Aut (B\otimes \mathcal{O}_{\infty })=\Aut (\mathcal{O}_{\infty }) \)
we get the desired entropy preserving covariant embedding. 
\end{proof}
\begin{subquestion} Let \( u\in A \) be a unitary and regard \( \Ad u \) as an automorphism of \( (A,\phi )*(B,\psi ) \).
As in Proposition \ref{prop:innerfreeeautomorphism}, is it true that \( \HT _{A*B}(\Ad u)=\HT (\Ad u|_{A}) \)?
Clearly, an affirmative answer to this question would imply that the automorphism in
Theorem \ref{theorem:cuntzalgebrahasallentropy} can be taken to be inner.\end{subquestion}

\subsection{Possible values of entropy of all automorphisms of a \protect\( C^{*}\protect \)--algebra.}

\begin{subdefn}
If \( A \) is an exact \( C^{*} \)--algebra then put
\[
TE(A)=\{t\in [0,\infty ]:{\textrm{there exists}}\, \alpha \in \Aut (A)\, {\textrm{such that}}\, \HT (\alpha )=t\}.\]
For $A$ unital, we can also consider the set
\[
TE_{Inn}(A)=\{t\in [0,\infty ]:{\textrm{there exists a unitary}}\, u\in A\, {\textrm{such that}}\, \HT ({\Ad }u)=t\}.\]
\end{subdefn}

The results of \S\ref{sec:htpresemb} imply, in particular, that \( TE(\mathcal{O}_{\infty })=[0,\infty ] \)
and \( TE_{Inn}(\mathcal{O}_{2})=[0,\infty ] \), since there are automorphisms of nuclear (in fact, abelian) algebras with any prescribed entropy. (See \cite{boyle-handleman:entropycantorset} for some nice examples.)  Clearly, if \( B \) contains
a projection, then \( TE(B\otimes A)\supset TE(A) \).
If $A$ and $B$ are both unital, we have \( TE_{Inn}(A\otimes B)\supset
TE_{Inn}(A)\cup TE_{Inn}(B) \). 

If \( B \) is any nuclear simple separable purely infinite \( C^{*} \)--algebra
then Kirchberg has shown that \( B\otimes \mathcal{O}_{\infty }\cong B \) (see
\cite{kirchberg-phillips:embeddingCuntz}).   Recently Kirchberg and R\o{}rdam
introduced a class of nonsimple purely infinite \( C^{*} \)--algebras (cf. \cite{kirchberg-rordam:nonsimplepureinf}).
 Though it is not yet known whether all of their nuclear purely infinite \( C^{*} \)--algebras
will absorb \( \mathcal{O}_{\infty } \), many examples are known. 

\begin{subthm}
Let \( B \) be any exact \( C^{*} \)--algebra which contains a projection.
Then \linebreak[4] \( TE(B\otimes \mathcal{O}_{\infty })=[0,\infty ] \). If \( B \) is
unital then \( TE_{Inn}(B\otimes \mathcal{O}_{2})=[0,\infty ] \). 

In particular, if \( B \) is nuclear, simple, separable and purely infinite
then \( TE(B)=[0,\infty ] \).
\end{subthm}

\bibliographystyle{amsplain}

\providecommand{\bysame}{\leavevmode\hbox to3em{\hrulefill}\thinspace}

\end{document}